\documentclass{amsart}
\usepackage{amsfonts}
\usepackage{amssymb}
\usepackage{amsmath}  %Overlap with newcommand \mod
\usepackage{verbatim}
\usepackage{stmaryrd}

\input xy
\xyoption{all}

\theoremstyle{plain} % style plain
\newtheorem{theorem}{Theorem}[section]
\newtheorem{lemma}[theorem]{Lemma}

\newtheorem{proposition}[theorem]{Proposition}
\newtheorem{corollary}[theorem]{Corollary}

\theoremstyle{definition} % style definition
\newtheorem{example}[theorem]{Example}
\newtheorem{definition}[theorem]{Definition}
\newtheorem{remark}[theorem]{Remark}
\newtheorem{remarks}[theorem]{Remarks}

\newcommand{\Ker}{\mbox{\rm Ker\,}}
\newcommand{\Coker}{\mbox{\rm Coker\,}}
\newcommand{\Ima}{\operatorname{\rm Im\,}}
\newcommand{\End}[1]{\operatorname{\rm End}_{#1}}
\newcommand{\Aut}[1]{\operatorname{\rm Aut}_{#1}}
\newcommand{\Hom}[1]{\operatorname{{\rm Hom}}_{#1}}
\newcommand{\HOM}[1]{\operatorname{{\rm HOM}}_{#1}}

\newcommand{\PHom}[1]{\operatorname{{\rm PHom}}_{#1}}

\newcommand{\Ext}[2]{\operatorname{\rm Ext}^{#1}_{#2}}

\newcommand{\dimv}{\underline{\dim}\,}

\newcommand{\add}{\mbox{{\rm add \!}}}
\newcommand{\codim}{\mbox{{\rm codim \!\!}}}
\newcommand{\GL}{\mbox{{\rm GL \!\!}}}
\newcommand{\MOD}{\mbox{{\rm mod \!}}}

\newcommand{\rep}{\mbox{{\rm rep \!\!}}}

\newcommand{\proj}{\operatorname{{\rm proj }}}
\newcommand{\perf}{\mbox{{\rm per \!}}}

\newcommand{\ind}[1]{\operatorname{\rm ind}_{#1}}

\newcommand{\Gr}[1]{\mbox{{\rm Gr}}_{#1}}

\newcommand{\K}{\operatorname{{\rm K}}}

\newcommand{\demo}[1]{\textsc{Proof} #1 \hfill $\Box$ \bigskip}

\newcommand{\cA}{\mathcal{A}}

\newcommand{\cC}{\mathcal{C}}
\newcommand{\cD}{\mathcal{D}}

\newcommand{\cF}{\mathcal{F}}

\newcommand{\cO}{\mathcal{O}}

\newcommand{\cU}{\mathcal{U}}
\newcommand{\cV}{\mathcal{V}}

\newcommand{\cX}{\mathcal{X}}
\newcommand{\cY}{\mathcal{Y}}
\newcommand{\cZ}{\mathcal{Z}}
\newcommand{\bC}{\mathbb{C}}
\newcommand{\bN}{\mathbb{N}}
\newcommand{\bP}{\mathbb{P}}
\newcommand{\bQ}{\mathbb{Q}}
\newcommand{\bZ}{\mathbb{Z}}
\newcommand{\bd}{\mathbf{d}}
\newcommand{\be}{\mathbf{e}}

\newcommand{\bi}{\mathbf{i}}
\newcommand{\bj}{\mathbf{j}}

\newcommand{\bm}{\mathbf{m}}

\begin{document}

\title[Generic bases for cluster algebras from the cluster category]{Generic bases for cluster algebras from the cluster category}
\author{Pierre-Guy Plamondon}
\email{pierre-guy.plamondon@unicaen.fr}
\address{ Laboratoire LMNO \\
          \'Equipe Alg\`ebre, G\'eom\'etrie et Logique \\
          Universit\'e de Caen \\
          F14032 Caen Cedex \\
          France 
   %Universit\'e Paris Diderot -- Paris 7\\
   %Institut de Math\'ematiques de Jussieu, UMR 7586 du CNRS\\
   %Case 7012\\
   %B\^atiment Chevaleret\\
   %75205 Paris Cedex 13\\
   %France 
   }

\thanks{The author was financially supported by a scholarship from the FQRNT}
%\date{\today}

\begin{abstract}
Inspired by recent work of Geiss--Leclerc--Schr\"oer, we use $\Hom{}$-finite cluster categories to give a good candidate set for a basis of (upper) cluster algebras with coefficients arising from quivers.  This set consists of generic values taken by the cluster character on objects having the same index.  If the matrix associated to the quiver is of full rank, then we prove that the elements in this set are linearly independent. If the cluster algebra arises from the setting of Geiss--Leclerc--Schr\"oer, then we obtain the basis found by these authors.  We show how our point of view agrees with the spirit of conjectures of Fock--Goncharov concerning the parametrization of a basis of the upper cluster algebra by points in the tropical $\cX$-variety.
\end{abstract}

\maketitle

\tableofcontents

%---------------------------------------------------------------------------------
\section{Introduction}
One of the main motivations of S.~Fomin and A.~Zelevinsky for introducing cluster algebras in \cite{FZ02} was the search for a combinatorial framework in which one could study the canonical bases of M.~Kashiwara \cite{Kashiwara90} and G.~Lusztig \cite{Lusztig90}.  Recent results of C.~Geiss, B.~Leclerc and J.~Schr\"oer \cite{GLS10}, who prove that coordinate rings of certain algebraic varieties have a natural cluster algebra structure and find a basis for them, give ample justification to this approach.  The problem of finding ``good'' bases for cluster algebras is thus central in the theory.  These bases should, as conjectured already in \cite{FZ02}, contain the cluster monomials.  Good bases for cluster algebras were also constructed by P.~Sherman and A.~Zelevinsky \cite{SZ04}, P.~Caldero and B.~Keller \cite{CK08}, G.~Dupont \cite{Dupont08} \cite{Dupont11}, D.~Hernandez and B.~Leclerc \cite{HL10}, H.~Nakajima \cite{Nak11}, M.~Ding, J.~Xiao and F.~Xu \cite{DXX09}, G.~Cerulli Irelli \cite{Cerulli09}, G.~Cerulli Irelli and F.~Esposito \cite{CE11}, G.~Dupont and H.~Thomas \cite{DT11} and G.~Musiker, R.~Schiffler and L.~Williams \cite{MSW11}. 

In their paper \cite{GLS10}, C.~Geiss, B.~Leclerc and J.~Schr\"oer find bases for a certain class of (upper) cluster algebras and provide a candidate for a basis in general. In this paper, inspired by their ideas, we use cluster categories to give another realization of this candidate set.  We prove that its elements are linearly independent when the defining matrix is of full rank, and that it coincides with the basis of \cite{GLS10} when the cluster algebra arises from the setting studied therein.

We work in the setting of cluster categories.  Our results apply to  cluster algebras $\cA_{Q,F}$ associated to ice quivers $(Q,F)$ \cite{FK09} on which there exists a non-degenerate potential $W$ (in the sense of \cite{DWZ08}) making $(Q,W)$ Jacobi-finite.  In that setting, C.~Amiot's (generalized) cluster category \cite{Amiot08} is $\Hom{}$-finite and is known to categorify the associated cluster algebra (see, for instance, \cite{Amiot08} and \cite{Palu08}).  The cluster category contains a distinguished object $\Gamma$ which possesses nice properties (for instance, it is \emph{cluster-tilting}).

In their paper \cite{FG09}, V.~Fock and A.~Goncharov give a set of conjectures related to the problem of finding a basis for (upper) cluster algebras.  In particular, they conjecture that the set $\cX(\bZ^t)$ of tropical $\bZ$-points of the $\cX$-variety parametrizes a basis of the upper cluster algebra $\cA^+_{Q,F}$.  There are natural isomorphisms
\begin{displaymath}
	\cX(\bZ^t) \textrm{ (for the opposite quiver) } \cong \bZ^n \cong K_0(\add\Gamma)
\end{displaymath}
which link the point of view of $\cX$-varieties to that of cluster categories; these isomorphisms commute with mutation on both sides.

 Each object $M$ of the cluster category has an \emph{index} $\ind{\Gamma}M$ (as defined in \cite{DK08}) which is an element of $K_0(\add \Gamma)$.  Our main theorem states that a good candidate for a basis of $\cA_{Q,F}$ is parametrized by the set of indices via the cluster character of Y.~Palu \cite{Palu08}.   Let $\bP$ be the tropical semifield generated by the frozen variables.

\begin{theorem}\label{theo::1}
There exists a canonical map
\begin{displaymath}
	I : K_0(\add \overline{\Gamma}) \longrightarrow \cA^+_{Q,F},
\end{displaymath}
where $\overline{\Gamma}$ is defined in section \ref{sect::reduction}. If the matrix of $(Q,F)$ is of full rank, then the elements in the image of $I$ are linearly independent over $\bZ\bP$.  If $(Q,W)$ arises from the setting of \cite{GLS10}, then the image of $I$ is the basis of the cluster algebra $\cA_Q$ found in that paper.
\end{theorem}
When $Q$ has no frozen vertices, the map $I$ sends an element $[T_0] - [T_1]$ to the generic value taken by the cluster character on cones of morphisms in $\Hom{\cC}(T_1, T_0)$.  It was considered by G.~Dupont in \cite{Dupont11}.

In their construction of a basis for cluster algebras, the authors of \cite{GLS10} consider \emph{strongly reduced components} of the variety $\rep(A)$ of finite-dimensional representations of some finite-dimensional algebra $A$.  It so happens that we can recover all such components from the set of indices $K_0(\proj A)$.

\begin{theorem}\label{theo::2}
Let $A$ be a finite-dimensional algebra given by a quiver with relations.  Then there exists a canonical surjection
\begin{displaymath}
	\Psi : K_0(\proj A) \longrightarrow \{\textrm{strongly reduced components of } \rep(A) \}.
\end{displaymath}
Two elements $\delta$ and $\delta'$ have the same image by $\Psi$ if, and only if, their canonical decompositions (in the sense of H.~Derksen and J.~Fei \cite{DF09}, see section \ref{sect::canonical decomposition}) can be written as
\begin{displaymath}
	\delta = \delta_1 \oplus \overline{\delta} \quad \textrm{and} \quad \delta' = \delta'_1 \oplus \overline{\delta},
\end{displaymath}
with $\delta_1$ and $\delta'_1$ non-negative.
\end{theorem}
Note that in the setting of the theorem, $K_0(\proj A)$ is isomorphic to $K_0(\add \Gamma)$, so the notation is coherent with that of Theorem \ref{theo::1}.

Elements of the cluster algebra can be mutated using the rules defined by S.~Fomin and A.~Zelevinsky \cite{FZ02}.  Elements of $K_0(\add \Gamma)$ can also be mutated (in a way which we will make precise).  The map $I$ of Theorem \ref{theo::1} is well-behaved with respect to those  different mutations, as conjectured in \cite[Conjecture 9.2]{Dupont11}.

\begin{theorem}\label{theo::3}
The map $I$ commutes with mutation, and thus its image is contained in the upper cluster algebra.
\end{theorem}

The paper is organized as follows.  We first prove \ref{theo::2} in section \ref{sect::}.  Then we define the map $I$ in section \ref{sect::generic value}, and then prove Theorem \ref{theo::3}.  Section \ref{sect::1} is devoted to the proof of Theorem \ref{theo::1}.  We end the paper with an example of a $\Hom{}$-finite cluster category for which the image of the map $I$ is not contained in the cluster algebra, and in which there are cluster-tilting objects that are not related by a sequence of mutations.

%-----------------------------------------------------------
\section{Projective presentations and strongly reduced components}\label{sect::}

\subsection{Varieties of representations}
Let $Q = (Q_0, Q_1, s,t)$ be a finite quiver, that is, an oriented graph with finitely many vertices and arrows.  We denote by $kQ$ its path algebra.  Let $I$ be an admissible ideal of $kQ$, and let $A=kQ/I$ be a finite-dimensional algebra (for general background on quivers and path algebras, we refer the reader to the book \cite{ASS}).

Let $\bd$ be a dimension vector for $Q$, that is, an element of $\bN^{Q_0}$.  The variety $\rep_{\bd}(A)$ is the affine variety whose points are representations of $Q$ with underlying space $\prod_{i\in Q_0} k^{d_i}$ satisfying the relations in $I$ ; it is realized as a Zariski-closed subset of the affine space $\prod_{a\in Q_1} \Hom{k}(k^{d_{s(a)}}, k^{d_{t(a)}})$.

We denote by $\rep(A)$ the disjoint union of all $\rep_{\bd}(A)$ as $\bd$ takes all possible values in $\bN^{Q_0}$.  For general background on varieties of representations, we refer the reader to \cite{C93}.

The algebraic group $\GL_{\bd}$ is defined to be $\prod_{i\in Q_0} \GL_{d_i}$.  It acts on $\rep_{\bd}(A)$ thus: for any $(g_i) \in \GL_{\bd}$ and any $(\varphi_a) \in \rep_{\bd}(A)$, $(g_i) (\varphi_a) = (g_{t(a)}\varphi_a(g_{s(a)})^{-1})$.  The orbit of a representation $M$ under the action of $GL_{\bd}$ is the set of representations with underlying space $\prod_{i\in Q_0} k^{d_i}$ isomorphic to $M$.

We will need the following information on the dimension of morphism and extension spaces, and on minimal projective presentations.

\begin{lemma}[Lemma 4.2 of \cite{CS02}]\label{lemm::CS}
The functions
\begin{displaymath}
	\rep_{\bd_1}(A) \times \rep_{\bd_2}(A) \longrightarrow \bZ
\end{displaymath}
sending a pair $(M_1, M_2)$ to the dimensions of the vector spaces $\Hom{A}(M_1, M_2)$ and $\Ext{1}{A}(M_1, M_2)$ are upper semicontinuous.
\end{lemma}

\begin{corollary}\label{coro::presentations}
Let $\cZ$ be an irreducible component of $\rep_{\bd}(A)$.  There exist finitely generated projective $A$-modules $P_1$ and $P_0$ and a dense open subset $\cU$ of $\cZ$ such that any representation $M$ in $\cU$ admits a minimal projective presentation of the form
\begin{displaymath}
	P_1 \longrightarrow P_0 \longrightarrow M \longrightarrow 0.
\end{displaymath}
\end{corollary}
\demo{ Given any representation $M$ and a minimal projective presentation $P_1 \longrightarrow P_0 \longrightarrow M \longrightarrow 0$, the multiplicities of an indecomposable projective $Q$ in $P_0$ and $P_1$ are given by the dimensions of $\Hom{A}(M,S)$ and $\Ext{1}{A}(M,S)$, respectively, where $S$ is a simple module whose projective cover is $Q$.

Restrict the maps of Lemma \ref{lemm::CS} to $\cZ \times \{S\}$.  The restrictions are still upper semicontinuous.  Therefore the subsets of $\cZ$ on which these functions take their minimal values are (dense) open subsets of $\cZ$. Their intersection is a dense open subset of $\cZ$ on which the functions
\begin{displaymath}
	\dim \Hom{A}(?, S) \quad \textrm{and} \quad \dim \Ext{1}{A}(?, S)
\end{displaymath}
are constant.  This proves the result.
}

We now introduce a slight modification of a definition of \cite[Section 7.1]{GLS10}.  Let $\cZ$ be an irreducible component of $\rep(A)$.  There is an open dense subset $\cU$ of $\cZ$ and positive integers $h(\cZ)$, $h'(\cZ)$, $e(\cZ)$ and $c(\cZ)$ such that, for any $M$ in $\cU$,
\begin{enumerate}
	\item $\dim \Hom{A}(M, \tau M) = h(\cZ)$ and $\dim \Hom{A}(\tau^{-1}M, M) = h'(\cZ)$ ;
	\item $\dim \Ext{1}{A}(M,M) = e(\cZ)$ ; and
	\item $\codim_{\cZ}(\GL_{\bd} M) = c(\cZ)$,
\end{enumerate}
where $\tau$ is the Auslander--Reiten translation (see, for example, \cite[Chapter IV]{ASS}).
Moreover, we have that $c(\cZ) \leq e(\cZ) \leq h(\cZ)$ and $c(\cZ) \leq e(\cZ) \leq h'(\cZ)$.

\begin{definition}\label{defi::strred}
\begin{itemize}
  \item (Section 7.1 of \cite{GLS10}) An irreducible component $\cZ$ of $\rep(A)$ such that $c(\cZ) = h'(\cZ)$ is \emph{strongly reduced}.
  \item An irreducible component $\cZ$ of $\rep(A)$ such that $c(\cZ) = h(\cZ)$ is \emph{dually strongly reduced}.
\end{itemize}
\end{definition}

\begin{remark}\label{rema::dually strongly reduced}
\begin{enumerate}

\item We use the terminology \emph{dually strongly reduced} to emphasize the following fact.  Let $D = \Hom{k}(?, k)$ be the usual duality functor.  It induces a (non-canonical) isomorphism of varieties $\rep(A) \cong \rep(A^{op})$.  Then an irreducible component $\cZ$ of $\rep(A)$ is strongly reduced if, and only if, its image is a dually strongly reduced component of $\rep(A^{op})$.  Indeed, for any finite-dimensional $A$-module $M$, we have that $\Ext{1}{A}(M,M) \cong \Ext{1}{A^{op}}(DM, DM)$, that $\dim \GL_{\bd}M = \dim \GL_{\bd}DM$, and that
\begin{displaymath}
	\Hom{A^{op}}(DM, \tau DM) = \Hom{A^{op}}(DM, D Tr D M) = \Hom{A^{op}}(DM, D \tau^{-1} M) \cong \Hom{A}(\tau^{-1}M, M),
\end{displaymath}
where $Tr$ is the transpose (again, we refer the reader to \cite[Chapter IV]{ASS}).

\item In the case where $A$ is the Jacobian algebra of a quiver with potential (see \cite{DWZ08}), the two definitions coincide.  Indeed, we have equalities
\begin{eqnarray*}
 \dim \Hom{A}(M, \tau M) & = & E^{proj}(M) \ (\textrm{by \cite[Corollary 10.9]{DWZ09}}) \\
                         & = & E^{proj}(\tau^{-1}M) \ (\textrm{by \cite[Corollary 7.5]{DF09}}) \\
                         & = &  \dim \Hom{A}(\tau^{-1}M, M) \ (\textrm{by \cite[Corollary 10.9]{DWZ09}}). 
\end{eqnarray*}
This can also be seen by using the cluster category $\cC$ of the quiver with potential, as defined in \cite{Amiot08}.  This category has a canonical cluster-tilting object $\Gamma$, and the functor $F = \Hom{\cC}(\Gamma, ?)$ induces an equivalence $\Hom{\cC}(\Gamma, ?)/(\Sigma\Gamma)\rightarrow \MOD A$ (\cite[Proposition 2.1(c)]{KR07}), such that $F(\Sigma X) = \tau (FX)$ (\cite[Section 3.5]{KR07}).  Then we have
\begin{eqnarray*}
\dim \Hom{A}(M, \tau M) & = & \dim \Hom{\cC}(\overline{M}, \Sigma \overline{M}) - \dim (\Sigma \Gamma)(\overline{M}, \Sigma\overline{M}) \\
                        & = & \dim \Hom{\cC}(\Sigma^{-1}\overline{M}, \overline{M}) - \dim (\Gamma)(\Sigma^{-1}\overline{M}, \overline{M}) \\
                        & = & \dim \Hom{\cC}(\Sigma^{-1}\overline{M}, \overline{M}) - \dim (\Sigma\Gamma)(\Sigma^{-1}\overline{M}, \overline{M}) \\
                        & = & \dim \Hom{A}(\tau^{-1}M, M),
\end{eqnarray*}
where $\overline{M}$ is a preimage of $M$ by the functor $F$, and the second-to-last equality is a consequence of \cite[Lemma 10]{Palu08}.

\end{enumerate}
\end{remark}

\subsection{Decomposition of projective presentations}\label{sect::canonical decomposition}

Let $A$ be a finite-dimensional algebra, and let $P'_1$, $P'_0$, $P''_1$ and $P''_0$ be finitely generated projective $A$-modules.

\begin{definition}[Definition 3.1 of \cite{DF09}]\label{defi::DF Einv}
For any $f'$ in $\Hom{A}(P'_1, P'_0)$ and any $f''$ in $\Hom{A}(P''_1, P''_0)$, define the space $E(f', f'')$ as
\begin{displaymath}
	E(f', f'') = \Hom{\K^b(\proj A)}(\Sigma^{-1}f', f''),
\end{displaymath}
where $f'$ and $f''$ are viewed as complexes in $\K^b(\proj A)$.  Define $E(f')$ to be $E(f', f')$.
\end{definition}

\begin{lemma}\label{lemm::Etau}
If $\xymatrix{P'_1 \ar[r]^{f'} & P'_0 \ar[r] & M' \ar[r] & 0}$ is a projective presentation, then 
\begin{displaymath}
	\dim E(f', f'') \geq \dim \Hom{A}(M'', \tau M'),
\end{displaymath}
where $M''$ is the cokernel of $f''$.  Equality holds if the presentation is minimal.
\end{lemma}
\demo{  Applying the right exact functor $D\Hom{A}(?, A) = D(?)^t$ to the presentation, we get an exact sequence
\begin{displaymath}
	\xymatrix{ 0\ar[r] & \tau M'\!\oplus I\ar[r] & D(P'_1)^t\ar[r] & D(P'_0)^t\ar[r] & D(M')^t\ar[r] & 0,
	}
\end{displaymath}
where $I$ is a finite-dimensional injective $A$-module which vanishes if the presentation is minimal.
We use the fact that the morphism of functors $D\Hom{A}(X,?) \rightarrow \Hom{A}(?, DX^t)$ is an isomorphism whenever $X$ is projective, and we get a commutative diagram with exact rows and vertical isomorphisms
\begin{displaymath}
	\xymatrix{ & & D\Hom{A}(P'_1, M'')\ar[d]^{\cong}\ar[r]^{D\phi} & D\Hom{A}(P'_0, M'')\ar[d]^{\cong} \\
	 0\ar[r] & \Hom{A}(M'', \tau M'\!\oplus I)\ar[r] & \Hom{A}(M'', D(P'_1)^t)\ar[r] & \Hom{A}(M'', D(P'_0)^t).
	}
\end{displaymath}
Therefore $D\Hom{A}(M'', \tau M'\!\oplus I)$ is isomorphic to the cokernel of $\phi$, which is in turn isomorphic to $E(f', f'')$ by \cite[Lemma 3.4]{DF09}.  This proves the inequality.  If the presentation is minimal, then $I$ vanishes and the equality holds.

}

We will need a result on the decomposition of general projective presentations, which follows from the work of H.~Derksen and J.~Fei on the one hand, and from that of W.~Crawley-Boevey and J.~Schr\"oer on the other hand.

For any $\delta$ in $K_0(\proj A)$, let $\PHom{A}(\delta)$ be the space $\Hom{A}(P^{\delta}_-, P^{\delta}_+)$, where $\delta = [P^{\delta}_+] - [P^{\delta}_-]$, and $P^{\delta}_+$ and $P^{\delta}_-$ have no non-zero direct factors in common.  If $[P^{\delta}_-] = 0$, then $\delta$ is called \emph{non-negative}.  Define $e(\delta, \delta')$ to be the minimal value of $\dim E(f,f')$ for $f \in \PHom{A}(\delta)$ and $f' \in \PHom{A}(\delta')$.

The vector $\delta$ is \emph{indecomposable} if a general element of $\PHom{A}(\delta)$ is indecomposable.  Its \emph{canonical decomposition} is $\delta_1 \oplus \ldots \oplus \delta_s$ if a general element of $\PHom{A}(\delta)$ has the form $f_1 \oplus \ldots \oplus f_s$, with $f_i \in \PHom{A}(\delta_i)$  and each $\delta_i$ is indecomposable \cite[Definition 4.3]{DF09}.

\begin{theorem}[Derksen--Fei]\label{theo::decomposition}
Any $\delta \in K_0(\proj A)$ admits a canonical decomposition $\delta_1 \oplus \ldots \oplus \delta_s$, where  $\delta_1, \ldots, \delta_s \in K_0(\proj A)$ are unique up to reordering.  A decomposition $\delta = \delta_1 \oplus \ldots \oplus \delta_s$ is the canonical decomposition if, and only if, each $\delta_i$ is indecomposable and $e(\delta_i, \delta_j) = 0$ whenever $i\neq j$.
\end{theorem}
\demo{The second statement is exacly \cite[Theorem 4.4]{DF09}.  To prove the first one, let $\bd = \dimv P_-^{\delta}$ and $\be = \dimv P_+^{\delta}$.  Then the orbit of $P_-^{\delta}$ (or $P_+^{\delta}$) is a dense open subset of an irreducible component $C_-$ of $\rep_{\bd}(A)$ (or $C_+$ of $\rep_{\be}(A)$, respectively), since $P_-^{\delta}$ is projective and thus has no self-extensions (see \cite[Corollary 1.2]{Gabriel74}).  Let
\begin{displaymath}
\rep_{\bd, \be}(A\overrightarrow{A_2}) = \{ (L,M,f) \big| \ L \in \rep_{\bd}(A), M \in \rep_{\be}(A), f \in \Hom{A}(L,M) \}.
\end{displaymath}
Then we can view $\PHom{A}(\delta) = \Hom{A}(P^{\delta}_-, P^{\delta}_+)$ as an irreducible subvariety of the affine variety $\rep_{\bd, \be}(A\overrightarrow{A_2})$.  Let $C$ be an irreducible component of $\rep_{\bd, \be}(A\overrightarrow{A_2})$ which contains $\PHom{A}(\delta)$.  By \cite[Theorem 1.1]{CS02}, there is a dense open subset $\cU$ of $C$ and indecomposable irreducible components $C_1, \ldots, C_s$ of $\rep(A\overrightarrow{A_2})$ such that $\cU \subset C_1 \oplus \ldots \oplus C_s$.  We have a diagram
\begin{displaymath}
	\xymatrix{ & C\ar[dl]_{\pi_-}\ar[dr]^{\pi_+}  & \\
	           C_- &   & C_+
	}
\end{displaymath}
where $\pi_-$ and $\pi_+$ are the natural projections; their images intersect the orbits of $P_-^{\delta}$ and $P_+^{\delta}$, respectively.  Thus the preimages of these open orbits are dense open subsets of $C$, whose common intersection with $\cU$ is a dense open subset $\cV$ of $C$.  Now $\cV\cap \PHom{A}(\delta)$ is non-empty, and is thus dense and open in $\PHom{A}(\delta)$.  The inclusion $\cV \cap \PHom{A}(\delta) \subset C_1\oplus \ldots \oplus C_s$ induces the canonical decomposition $\delta_1 \oplus \ldots \oplus \delta_s$ of $\delta$.

}

\begin{corollary}\label{coro::minimal}
If the canonical decomposition of $\delta$ has no non-negative factors, then a general element in $\PHom{A}(\delta)$ is a minimal projective presentation.
\end{corollary}
	
We end this section by a useful result concerning ``orthogonal'' presentations.	
	
\begin{definition}
Two elements $\delta' = \sum_{i=1}^n g'_i [P_i]$ and $\delta'' = \sum_{i=1}^n g''_i [P_i]$ of $K_0(\proj A)$ are \emph{sign-coherent} if, for all $i$ from $1$ to $n$, the integers $g'_i$ and $g''_i$ are both non-negative or both non-positive.
\end{definition}

\begin{lemma}\label{lemm::signcoherence}
If $e(\delta', \delta'') = e(\delta'', \delta') =0$, then $\delta'$ and $\delta''$ are sign-coherent.
\end{lemma}	
\demo{ Assume that $\delta'$ and $\delta''$ are not sign-coherent.  Then, without loss of generality, we can write $\delta' = [P'_0] - [P'_1]$ and $\delta'' = [P''_0] - [P''_1]$ and assume that $P'_0$ and $P''_1$ have a direct factor in common.  We also assume that $P'_0$ and $P'_1$ share no direct factors, and that $P''_0$ and $P''_1$ share no direct factors.  Finally, note that any morphism between objects of $\proj A$ can be written as a matrix whose entries are elements of $A$.  Since $P'_0$ and $P''_1$ have a common direct factor, there exists a morphism $h$ from $P''_1$ to $P'_0$ whose matrix form contains entries not in the radical of $A$ (for instance, if $R$ is the common direct factor, take the composition $P''_1 \rightarrow R \rightarrow P'_0$, where the left morphism is a retraction and the right one is a section). 

Consider the diagram
\begin{displaymath}
	\xymatrix{  0\ar[r] & P''_1 \ar[r]^{f''}\ar[d]^{h}\ar[dl]_s & P''_0\ar[dl]^t \\
	            P'_1\ar[r]^{f'} & P''_0\ar[r] & 0
	}
\end{displaymath}
where $s$, $t$, $f'$ and $f''$ are any morphisms.  Since $P'_0$ and $P'_1$ share no direct factors, the matrix of $f'$ has all its entries in the radical of $A$.  The same holds for the matrix of $f''$.  Therefore the sum $f's + tf''$ has all its entries in the radical of $A$, and cannot be equal to $h$.  This proves that $E(f'', f') = \Hom{K^b(\proj A)}(\Sigma^{-1}f'', f')$ cannot vanish for any $f'$ and $f''$, and thus that $e(f'', f')$ is non-zero.
}

%............................................................
\subsection{Morphic cokernels}\label{sect::cokernels}

Let $A$ be a finite-dimensional $k$-algebra as before, and let $\bd$, $\bd_1$ and $\bd_0$ be dimension vectors.  Define the affine varieties
\begin{displaymath}
\rep_{\bd_1, \bd_0}(A\overrightarrow{A_2}) = \{ (L,M,f) \big| \ L \in \rep_{\bd_1}(A), M \in \rep_{\bd_0}(A), f \in \Hom{A}(L,M) \}
\end{displaymath}
\begin{displaymath}
\rep_{\bd_1, \bd_0}(A\overrightarrow{A_2})_{\bd} = \{ (L,M,f) \in \rep_{\bd_1, \bd_0}(A\overrightarrow{A_2}) \big| \ \dimv \Coker f = \bd \}.
\end{displaymath}

The latter is a locally closed subset of the former.  The symbol $\overrightarrow{A_2}$ stands for the quiver $1\rightarrow 2$ ; elements of the above sets are $A$-module-valued representations of $\overrightarrow{A_2}$.

Fix bases $\{ u_1, u_2, \ldots, u_{\ell} \}$ and $\{ v_1, v_2, \ldots, v_m \}$ of $\prod_{i\in Q_0} k^{d_{1, i}}$ and $\prod_{i\in Q_0} k^{d_{0, i}}$, respectively (these are the underlying vector spaces of $L$ and $M$) ; choose the basis vectors so that they all lie in some $k^{d_{e, i}}$ for $e = 0,1$ and $i\in Q_0$.  For any subset $\bi$ of $\{1, 2, \ldots, m\}$, let $N_{\bi}$ be the vector space generated by $\{v_i | \ i \in \bi \}$, and let
\begin{displaymath}
	E_{\bi} = \{ (L,M,f) \in \rep_{\bd_1, \bd_0}(A\overrightarrow{A_2}) \big| \ M \cong N_{\bi} \oplus \Ima f \ \textrm{as a vector space} \}.
\end{displaymath}

Notice that $\rep_{\bd_1, \bd_0}(A\overrightarrow{A_2})$ is the union of the $E_{\bi}$, and that each $E_{\bi}$ is contained in $\rep_{\bd_1, \bd_0}(A\overrightarrow{A_2})_{\bd}$ for some dimension vector $\bd$.  Notice also that $E_{\bi}$ is the intersection of an open subset with a closed subset.  Indeed, an element $(L, M, f)$ of $\rep_{\bd_1, \bd_0}(A\overrightarrow{A_2})$ lies in $E_{\bi}$ if and only if the following two conditions are satisfied (here we write $f$ as $(a_{ij})$ in matrix form with respect to the fixed bases) :

\begin{description}
	\item[(a)] There exists a subset $\bj$ of $\{1, \ldots, \dim P_1\}$ such that $|\bj| = m - |\bi|$ and the submatrix $(a_{ij})_{i\notin \bi, j \in \bj}$ has a non-zero determinant.  This condition defines an open subset.
	\item[(b)] For any $i_0 \in \bi$, and any subset $\bj$ of $\{ 1, \ldots, \dim P_1 \}$ such that $|\bj| = m - |\bi| + 1$, the submatrix $(a_{ij})$, where $j\in \bj$ and $i$ is either $i_0$ or not in $\bi$, has vanishing determinant.  This condition defines a closed subset.
\end{description}

In particular, if $E_{\bi}$ is contained in $\rep_{\bd_1, \bd_0}(A\overrightarrow{A_2})_{\bd}$, then it is open inside it, since the second condition is then automatically satisfied.
The next statement is a slight generalization of a result of \cite[Lemma 4]{Palu09}.

\begin{lemma}\label{lemm::morphism}
Assume that $E_{\bi}$ is contained in $\rep_{\bd_1, \bd_0}(A\overrightarrow{A_2})_{\bd}$.  Then there exists a morphism of varieties
  \begin{displaymath}
	  \Phi : E_{\bi} \longrightarrow \rep_{\bd}(A)
  \end{displaymath}
such that $\Phi(f)$ is isomorphic to $\Coker f$ for any element $f$ of $E_{\bi}$.
\end{lemma}
\demo{Let $(L, M, f)$ be an element of $E_{\bi}$.  We define $\Phi(L,M,f)$, as a vector space, to be the quotien of $M$ by $\Ima f$, that is, $N_{\bi}$. Let us define the $A$-module structure.

Let $\Omega_{\bj}$ be the open subset of $E_{\bi}$ consisting of maps satisfying condition (a) above for some fixed $\bj$. Then the $\Omega_{\bj}$ form an open cover of $E_{\bi}$.  

We define $\Phi$ on $\Omega_{\bj}$ as follows.  Assume that $(L, M, f)$ lies in $\Omega_{\bj}$.  Let $D = (a_{ij})_{i\notin \bi, j\in \bj}$ and $C = (a_{ij})_{i\in \bi, j\in \bj}$ ; then $D$ is invertible by condition (a).  

Let $b$ be an element of $A$.  We will define the matrix of the action of $b$ on $\Phi(L,M,f)$.  We do this through the following diagram :
\begin{displaymath}
	\xymatrix{ E_{\bi} \ar[r]^{I_{\bi}} & M \ar[r]^{\rho_M(b)} & M\ar[r]^{P_{\bi} \phantom{xx}} &  E_{\bi} \oplus E_{\bi'} \ar[r]^{F} & E_{\bi} \oplus \Ima f \ar[r]^{\phantom{xxx} \pi} & E_{\bi}.
	}
\end{displaymath}
Here $I_{\bi}$ is the natural inclusion; $\rho_M(b)$ is the action of $b$ on $M$; $P_{\bi}$ is the permutation matrix putting the basis vectors $v_i$, $i\in \bi$, before the others;  $F$ is a base change matrix given by
\begin{displaymath}
	\left(
	\begin{array}{cc}
1 & -CD^{-1} \\
0 & D^{-1}
\end{array} \right)
\end{displaymath}
and $\pi$ is the natural projection, given by $(1, 0)$.  Thus the action of $b$ on $\Phi(L,M,f)$ is given by the matrix
\begin{displaymath}
	(1, -CD^{-1})P_{\bi}(\rho_M(b))I_{\bi}.
\end{displaymath}

We have thus defined the action of $\Phi$ on $\Omega_{\bj}$.  This definition does not depend on $\bj$ ; indeed, assume that $(L,M,f)$ is also in $\Omega_{\bj'}$.  Let $D' = (a_{ij})_{i\notin \bi, j\in \bj'}$ and $C' = (a_{ij})_{i\in \bi, j\in \bj'}$.  Then $CD^{-1} = C'(D')^{-1}$.  To see this, notice that condition (b) above implies that any line of the matrix $(a_{ij})$ which is in $\bi$ is a linear combination of the ones not in $\bi$.  Therefore there exists a matrix $K$ such that $(a_{ij})_{i\in \bi} = K(a_{ij})_{i\notin \bi}$.  Therefore $C = KD$ and $C' = KD'$, and we get the desired equality.

Therefore $\Phi$ is well-defined on an open cover of $E_{\bi}$, and it is thus a morphism of varieties.

}

%....................
\subsection{Codimensions of orbits}

In the preceding section we have defined a morphism $\Phi : E_{\bi} \longrightarrow \rep_{\bd}(A)$.  Recall that $\Phi(L,M,f)$ is isomorphic to $\Coker f$, and that  an open cover of $\rep_{\bd_1, \bd_0}(A\overrightarrow{A_2})_{\bd}$ is formed by such $E_{\bi}$'s.

Define $\GL_{\bd_1,\bd_0}$ as the algebraic group $\GL_{\bd_1} \times \GL_{\bd_0}$.  Then the group $\GL_{\bd_1,\bd_0}$ acts on $\rep_{\bd_1, \bd_0}(A\overrightarrow{A_2})_{\bd}$ thus : for any $(g_1, g_0) \in \GL_{\bd_1,\bd_0}$ and any $(L,M,f) \in \rep_{\bd_1, \bd_0}(A\overrightarrow{A_2})_{\bd}$, we have that $(g_1, g_0)(L,M,f) = (g_1 L, g_0 M, g_0 f g_1^{-1})$.

\begin{lemma}\label{lemm::dominant2}
Let $(L,M,f)$ be an element of $E_{\bi}$.  Then the orbit of $\Phi(L,M,f)$ in $\rep_{\bd}(A)$ is equal to the image by $\Phi$ of the intersection of the orbit of $(L,M,f)$ with $E_{\bi}$.  In short,
\begin{displaymath}
	\cO_{\Phi(L,M,f)} = \Phi(E_{\bi} \cap \cO_{(L,M,f)}).
\end{displaymath}
\end{lemma}
\demo{  Let $b$ be an element of $A$.  Then we showed in lemma \ref{lemm::morphism} that $b$ acts on $\Phi(L,M,f)$ by the matrix
\begin{displaymath}
	(1, -CD^{-1})P_{\bi}(\rho_M(b))I_{\bi}.
\end{displaymath}
Let $\gamma$ be an element of $\GL_{\bd}$.  Then the action of $b$ on $\gamma \Phi(L,M,f)$ is
\begin{displaymath}
	\gamma (1, -CD^{-1})P_{\bi}(\rho_M(b))I_{\bi} \gamma^{-1}.
\end{displaymath}

Consider the element $G = (1, \overline{\gamma})$ of $\GL_{\bd_1, \bd_0}$, where 
\begin{displaymath}
	G = P_{\bi}^{-1}\left(
	\begin{array}{cc}
\gamma & 0 \\
0 & 1
\end{array} \right)P_{\bi}.
\end{displaymath}

Then $G(L,M,f) = (L, \overline{\gamma}M, \overline{\gamma}f)$.  In matrix form, we have that
\begin{displaymath}
	\overline{\gamma} f = P_{\bi}^{-1}\left(
	\begin{array}{cc}
\gamma & 0 \\
0 & 1
\end{array} \right)P_{\bi} P_{\bi}^{-1} \left(
	\begin{array}{cc}
C & A \\
D & B
\end{array} \right) P_{\bj} = P_{\bi}^{-1}\left(
	\begin{array}{cc}
\gamma C & \gamma A \\
 D &  B
\end{array} \right)P_{\bj}.
\end{displaymath}

Therefore $G(L,M,f)$ is still in $E_{\bi}$, and $b$ acts on $\Phi(G(L,M,f))$ by
\begin{eqnarray*}
	& & (1, -\gamma CD^{-1})P_{\bi}(\rho_{\overline{\gamma}M}(b))I_{\bi} \\	
	&=& (1, -\gamma CD^{-1})P_{\bi}P_{\bi}^{-1}\left(
	       \begin{array}{cc}
          \gamma & 0 \\
          0 & 1
          \end{array} \right)P_{\bi}(\rho_{M}(b)) P_{\bi}^{-1} \left(
	        \begin{array}{cc}
          \gamma^{-1} & 0 \\
          0 & 1
          \end{array} \right)P_{\bi}I_{\bi} \\
  &=& (\gamma , -\gamma CD^{-1})P_{\bi}(\rho_{M}(b)) P_{\bi}^{-1} \left(
	        \begin{array}{cc}
          \gamma^{-1} & 0 \\
          0 & 1
          \end{array} \right) \left(
	        \begin{array}{c}
          1 \\
          0
          \end{array} \right) \\
  &=&     \gamma(1 , - CD^{-1})P_{\bi}(\rho_{M}(b)) P_{\bi}^{-1} \left(
	        \begin{array}{c}
          1 \\
          0
          \end{array} \right)\gamma^{-1} \\
  &=&     \gamma(1 , - CD^{-1})P_{\bi}(\rho_{M}(b)) I_{\bi}\gamma^{-1}.
\end{eqnarray*}

Therefore $\Phi(G(L,M,f)) = \gamma\Phi(L,M,f)$.  This proves that we have an inclusion $\cO_{\Phi(L,M,f)} \subset \Phi(E_{\bi} \cap \cO_{(L,M,f)})$.

The other inclusion follows from the fact that if $(L',M',f')$ lies in the orbit of $(L,M,f)$, then the cokernels of $f$ and $f'$ are isomorphic.  This proves the lemma.

}

In the course of this section, we will be relying heavily on the following theorem on dimensions, borrowed from the book \cite{Borel}.

\begin{theorem}[Theorem AG.10.1 of \cite{Borel}]\label{theo::dim}
Let $a:X\rightarrow Y$ be a dominant morphism of irreducible varieties.  Let $W$ be an irreducible closed subvariety of $Y$ and let $Z$ be an irreducible component of $a^{-1}(W)$.

There exists an open dense subset $U$ of $Y$ (depending only on $a$) such that 
\begin{itemize}
	\item $U \subset a(X)$, and
	\item if $Z$ and $a^{-1}(U)$ have non-empty intersection, then $\codim_X Z = \codim_Y W$.
\end{itemize}
\end{theorem}

For the next lemma, we shall make the following identifications and definitions:  
\begin{eqnarray*}
	\Hom{A}(L_0,M_0) &=& \{ (L,M,f) \in \rep_{\bd_1, \bd_0}(A\overrightarrow{A_2}) \big| \ L = L_0, M = M_0 \};\\
	\Hom{A}(L_0,M_0)_{\bd} &=& \Hom{A}(L_0,M_0) \cap \rep_{\bd_1, \bd_0}(A\overrightarrow{A_2})_{\bd};\\
	\HOM{A}(L_0,M_0) &=& \{ (L,M,f) \in \rep_{\bd_1, \bd_0}(A\overrightarrow{A_2}) \big| \ L \cong L_0, M \cong M_0 \};\\
	\HOM{A}(L_0,M_0)_{\bd} &=& \HOM{A}(L_0,M_0) \cap \rep_{\bd_1, \bd_0}(A\overrightarrow{A_2})_{\bd};\\
	\GL_{L_0,M_0} &=& \Aut{A}(L_0) \times \Aut{A}(M_0).
\end{eqnarray*}
We shall denote an element $(L_0, M_0,f)$ of $\Hom{A}(L_0, M_0)$ simply by the morphism $f$.  Note that the  first and the third varieties are irreducible; indeed, the first one is a vector space, and the third one is $GL_{\bd_1, \bd_0}\Hom{A}(L_0,M_0)$, which is irreducible.  Note that $\GL_{L_0,M_0}$ acts on $\Hom{A}(L_0,M_0)$.

Notice that, inside $\Hom{A}(L,M)$ and $\HOM{A}(L,M)$, the subsets of the $(L,M,f)$ such that $f$ is of maximal rank are open subsets, and the cokernels of those $f$ all have the same dimension vector.  We denote those  subsets by $\Hom{A}(L,M)_{max}$ and $\HOM{A}(L,M)_{max}$.

\begin{lemma}\label{lemm::orbits}
  Fix $(L_0, M_0) \in \rep_{\bd_1}(A) \times \rep_{\bd_0}(A)$.  Let $\bi$ be such that $E_{\bi}$ intersects $\HOM{A}(L_0, M_0)_{\bd}$, and consider the morphism $\Phi : E_{\bi} \rightarrow \rep_{\bd}(A)$ defined above.  There exists an open subset $\cV$ of $E_{\bi}\cap\HOM{A}(L_0, M_0)_{\bd} $ such that for any $(L,M,f)$ in $\cV$, the following properties hold.
\begin{enumerate}
	\item If $\cF$ is an irreducible component of $\HOM{A}(L,M)_{\bd}$ which contains $\cO_{(L,M,f)}$, then 
	  \begin{displaymath}
	    \codim_{\cF}\cO_{(L,M,f)} =\codim_{\Phi(E_{\bi}\cap \cF)}\cO_{\Phi(L,M,f)}.
    \end{displaymath}
    In particular, if $\HOM{A}(L,M)_{\bd} = \HOM{A}(L,M)_{max}$, then
    \begin{displaymath}
	    \codim_{\cY}\cO_{(L,M,f)} =\codim_{\Phi(E_{\bi}\cap \cY)}\cO_{\Phi(L,M,f)},
    \end{displaymath}
    where $\cY = \HOM{A}(L,M)$.
  
  \item With the same notation as in (1), and letting $\cX = \Hom{A}(L,M)$, we have that
    \begin{displaymath}
	    \codim_{\cX}\cO_f = \codim_{\cY}\cO_{(L,M,f)}.
    \end{displaymath}
\end{enumerate}
\end{lemma}
\demo{
We first prove (1).  Consider the following commuting diagram :
\begin{displaymath}
	\xymatrix{ \overline{\cO_{(L,M,f)}} \ar@{_{(}->}[d] & \overline{E_{\bi} \cap \cO_{(L,M,f)}}\ar@{_{(}->}[d] \ar[l] \ar[r]^{\Phi} & \overline{\cO_{\Phi(L,M,f)}}\ar@{_{(}->}[d] \\
	              \cF       & E_{\bi} \cap \cF \ar@{_{(}->}[l]\ar[r]^{\Phi}  &  \Phi(E_{\bi} \cap \cF). 
	}
\end{displaymath}

The three varieties in the lower row are irreducible.  Since $E_{\bi}\cap\cF$ is a dense open subset of $\cF$, the lower-left morphism is dominant.  The lower-right morphism is also dominant (see, for instance, \cite[AG.10.2]{Borel}).  So we can apply the dimension theorem \ref{theo::dim}; if $\cU_1\subset \cF$ and $\cU_2 \subset \Phi(E_{\bi} \cap \cF)$ are the open subsets described by the theorem, let $\cV$ be the intersection of their preimages in $E_{\bi} \cap \cF$.  

Now, $\overline{E_{\bi} \cap \cO_{(L,M,f)}}$ is an irreducible component of $\Phi^{-1}(\overline{\cO_{\Phi(L,M,f)}})$ thanks to Lemma \ref{lemm::dominant2}.  Moreover, $\overline{E_{\bi} \cap \cO_{(L,M,f)}}$ is the preimage of $\overline{\cO_{(L,M,f)}}$ by the inclusion.  Thus we can apply the dimension theorem and get
\begin{eqnarray*}
	\codim_{\cF}\cO_{(L,M,f)} &=& \codim_{E_{\bi} \cap \cF}E_{\bi} \cap \cO_{(L,M,f)} \\
	  &=&  \codim_{\Phi(E_{\bi} \cap \cF)}\cO_{\Phi(L,M,f)}.
\end{eqnarray*}
This proves the first result.

Let us now prove (2).  Consider the diagram
\begin{displaymath}
	\xymatrix{ \cO_f \ar@{^{(}->}[r] \ar@{_{(}->}[d]  &  \Hom{A}(L,M) \ar[r] \ar@{_{(}->}[d] & \{(L, M)\} \ar@{_{(}->}[d]\\
	           \cO_{(L,M,f)} \ar@{^{(}->}[r] & \HOM{A}(L,M) \ar[r]  & \cO_L \times \cO_M.
	}
\end{displaymath}

Two applications of the dimension theorem \ref{theo::dim} yields equalities
\begin{eqnarray*}
 \codim_{\cO_{(L,M,f)}} \cO_f &=& \dim \cO_L \times \cO_M \\
                              &=& \codim_\cY \cX,
\end{eqnarray*}
which in turn yields
\begin{displaymath}
	\codim_{\cX}\cO_f = \codim_{\cY}\cO_{(L,M,f)}.
\end{displaymath}

}

For any $f\in \Hom{A}(L, M)$, the action of $\GL_{L, M}$ induces a morphism
\begin{displaymath}
	\begin{array}{cccc}
	\pi: & \GL_{L, M} & \longrightarrow & \cO_f\\
	  & (g_1, g_0) & \longmapsto & g_0f(g_1)^{-1}
	\end{array}
\end{displaymath} 

which, in turn, induces a linear map on tangent spaces
\begin{displaymath}
	\begin{array}{cccc}
	d\pi: & \End{A}L \oplus \End{A}M & \longrightarrow & T_f(\cO_f)\\
	  & (h_1, h_0) & \longmapsto & fh_1 - h_0f.
	\end{array}
\end{displaymath} 
Here we view $T_f(\cO_f)$ as a subspace of $T_f(\Hom{A}(L, M))$, which we identify with the space $\Hom{A}(L, M)$.

\begin{lemma}\label{lemm::differential}
The map $d\pi$ is surjective.
\end{lemma}
\demo{The morphism $\pi$ is surjective by definition.  In particular, it is dominant.  Since we work over a field of characteristic zero, it is automatically separable.  It then follows from \cite[Proposition II.6.7 and AG.17.3]{Borel} that $d\pi$ is surjective.
}

%...............................
\subsection{Orbits and the E-invariant}
As before, let $A$ be a finite-dimensional $k$-algebra and let $P_1$ and $P_0$ be finitely generated $A$-modules.

\begin{lemma}\label{lemm::einj}
Let $f$ be any element of $\Hom{A}(P_1, P_0)$.  We have the equality
\begin{displaymath}
	\codim_{\cX} \cO_f = \dim E(f),
\end{displaymath}
where $\cX$ stands for $\Hom{A}(P_1, P_0)$.
\end{lemma}
\demo{ Consider the linear map
\begin{eqnarray*}
\psi : \Hom{A}(P_1, P_0) & \longrightarrow & E(f) \\
       g & \longmapsto & \overline{g},
\end{eqnarray*}
where $\overline{g}$ is the map from $\Sigma^{-1}f$ to $f$ in $\K^b(\proj A)$ given by
\begin{displaymath}
	\xymatrix{ \cdots\ar[r] & 0 \ar[r]\ar[d]^0 & P_1\ar[r]^{f}\ar[d]^g & P_0\ar[r]\ar[d]^0 & \cdots \\
	           \cdots\ar[r] & P_1\ar[r]^{f} & P_0\ar[r] & 0\ar[r] & \cdots.
	}
\end{displaymath}
The map $\psi$ is obviously surjective.  Moreover, its kernel is exactly
\begin{displaymath}
	T_f(\cO_f) = \{h_0f + fh_1 | h_i \in \End{A}(P_i) \},
\end{displaymath}
since this is the very definition of null-homotopic maps from $\Sigma^{-1}f$ to $f$ (the above equality follows from Lemma \ref{lemm::differential}).  Therefore
\begin{eqnarray*}
\codim_{\cX}\cO_f & = & \dim\big( \Hom{A}(P_1, P_0)/T_f(\cO_f) \big) \\
                  & = & \dim E(f).
\end{eqnarray*}  
}

%................................
\subsection{Proof of Theorem \ref{theo::2}}\label{sect::proof of 2}
We define the map
\begin{displaymath}
	\Psi : K_0(\add A) \longrightarrow \{\textrm{strongly reduced components of } \rep(A) \}.
\end{displaymath}
To do so, we first define a map
\begin{displaymath}
	\Psi' : K_0(\add A^{op}) \longrightarrow \{\textrm{dually strongly reduced components of } \rep(A^{op}) \}.
\end{displaymath}
For any element of the form $\delta = [P_0] - [P_1]$, where $P_1$ and $P_0$ are two projective modules over $A^{op}$ which share no non-zero direct factors, consider the morphism of varieties
\begin{displaymath}
	\Phi : E_{\bi}\cap \HOM{A^{op}}(P_0, P_1)_{max} \longrightarrow \rep_{\bd}(A^{op})
\end{displaymath}
constructed in section \ref{sect::cokernels}. 

By Lemma \ref{lemm::orbits}, there is a dense open subset of the set $\HOM{A^{op}}(P_0, P_1)_{max}$ such that, for any $(L,M,f)$ in that open subset, 
\begin{displaymath}
	\codim_{\cX}\cO_f = \codim_{\Phi(E_{\bi} \cap \cY)}\cO_{\Phi(L,M,f)}.
\end{displaymath} 
Now, by Lemma \ref{lemm::einj}, we have that
\begin{displaymath}
	\codim_{\cX} \cO_f = \dim E(f).
\end{displaymath}
Therefore, 
\begin{eqnarray*}
	\codim_{\cZ}\cO_{\Phi(L,M,f)} & \geq & \codim_{\Phi(E_{\bi} \cap \cY)} \cO_{\Phi(L,M,f)} \\
	                          & = & \dim E(f) \\
	                          & \geq & \dim \Hom{A^{op}}(\Phi(L,M,f), \tau \Phi(L,M,f)) \\
	                          & \geq & \codim_{\cZ}\cO_{\Phi(L,M,f)},
\end{eqnarray*}
where the third inequality follows from Lemma \ref{lemm::Etau}.  This implies that we have $\cZ = \overline{\Phi(E_{\bi} \cap \cY)}$, and that $\cZ$ is a dually strongly reduced component of $\rep(A^{op})$.

Define $\Psi'(\delta)$ to be this $\cZ$.

Now recall from Remark \ref{rema::dually strongly reduced} that the duality $D:\MOD(A) \rightarrow \MOD(A^{op})$ induces an isomorphism of varieties $\rep(A) \rightarrow \rep(A^{op})$ which sends strongly reduced components to dually strongly reduced ones; thus there is a strongly reduced component $\cZ_0$ of $\rep(A)$ corresponding to $\cZ$.  Moreover, $K_0(\add A) \cong K_0(\add A^{op})$ in a natural way; $\delta$ thus corresponds to some $\delta_0 \in K_0(\add A)$.  We define $\Psi(\delta_0)$ to be the strongly reduced component $\cZ_0$.

From this definition, it follows immediately that two elements $\delta$ and $\delta'$ of $K_0(\proj A^{op})$ have the same image by $\Psi'$ if, and only if, their canonical decompositions can be written as 
\begin{displaymath}
	\delta = \delta_1 \oplus \overline{\delta}, \quad \delta' = \delta'_1 \oplus \overline{\delta},
\end{displaymath}
with $\delta_1, \delta'_1$ non-negative, for non-negative factors do not affect the cokernels.

Let us now prove that $\Psi$ is surjective.  It suffices to show that $\Psi'$ is surjective.

Let $\cZ'$ be a dually strongly reduced component of $\rep(A^{op})$.  By Corollary \ref{coro::presentations}, there is a dense open subset $\cU$ of $\cZ'$ and there are finitely generated projective modules $P_1$ and $P_0$ such that every representation $M$ in $\cU$ admits a minimal projective presentation
\begin{displaymath}
	P_0 \longrightarrow P_1 \longrightarrow M \longrightarrow 0.
\end{displaymath}

Consider the locally closed subset $\HOM{A^{op}}(P_0, P_1)_{\bd}$.  There exists an irreducible component $\cF$  of it and an $\bi$ such that $\Phi(E_{\bi} \cap \cF) \cap \cU$ is dense in $\cZ'$. We get

\begin{eqnarray*}
\codim_{\cF}\cO_{(P_1, P_0, f)} & = & \codim_{\Phi(E_{\bi} \cap \cF)}\cO_{\Phi(P_0, P_1, f)} \textrm{ (Lemma \ref{lemm::orbits})} \\
 & = & \codim_{\cZ'}\cO_{\Phi(P_0, P_1, f)} \\
 & = & \dim \Hom{A^{op}}(\Phi(P_0, P_1, f), \tau \Phi(P_0, P_1, f)) \textrm{ ($\cZ'$ str. reduced)} \\
 & = & \dim E(f) \textrm{ (Lemma \ref{lemm::Etau})} \\
 & = & \codim_{\cX} \cO_f \textrm{ (Lemma \ref{lemm::einj})} \\
 & = & \codim_{\cY} \cO_{(P_0, P_1, f)} \textrm{ (Lemma \ref{lemm::orbits})}.
\end{eqnarray*}

Therefore $\cF$ is of codimension zero in $\HOM{A^{op}}(P_0, P_1)$.  Thus we have the equality $\HOM{A^{op}}(P_0, P_1)_{\bd} = \HOM{A^{op}}(P_0, P_1)_{max}$, and thus $\Psi'([P_0] - [P_1]) = \cZ'$.  This proves the surjectivity of the map $\Psi'$.

We would like to list now some properties that have been shown in the course of the proof.

\begin{corollary}
Let $A = kQ/I$ be a finite-dimensional $k$-algebra, and let $\cZ$ be an irreducible component of $\rep(A)$.
\begin{enumerate}
	\item If $\cZ$ is dually strongly reduced, then there exists a dense open subset $\cU$ of $\cZ$ and projective $A$-modules $P_1$ and $P_0$ (uniquely determined by $\cZ$) such that any $M$ in $\cU$ admits a minimal projective presentation
	\begin{displaymath}
	P_1 \longrightarrow P_0 \longrightarrow M \longrightarrow 0
\end{displaymath} 
and such that $P_1$ and $P_0$ have no direct factor in common.  Moreover, if $\cZ'$ is a different dually strongly reduced component, then the projectives $P'_1$ and $P'_0$ that it determines are such that $[P_0] - [P_1] \neq [P'_0] - [P'_1]$ in $K_0(\add A)$.

  \item If $\cZ$ is strongly reduced, then there exists a dense open subset $\cU$ of $\cZ$ and injective $A$-modules $I_0$ and $I_1$ (uniquely determined by $\cZ$) such that any $M$ in $\cU$ admits a minimal injective presentation
	\begin{displaymath}
	0 \longrightarrow M \longrightarrow I_0 \longrightarrow I_1 
\end{displaymath} 
and such that $I_0$ and $I_1$ have no direct factor in common.  Moreover, if $\cZ'$ is a different strongly reduced component, then the injectives $I'_0$ and $I'_1$ that it determines are such that $[I_0] - [I_1] \neq [I'_0] - [I'_1]$ in $K_0(\add DA^{op})$.
\end{enumerate}
\end{corollary} 
\demo{ Statement (2) is obtained from statement (1) by duality.  In statement (1), the existence of $\cU$ and of $P_1$ and $P_0$ follows from Corollary \ref{coro::presentations}. Following the end of the above proof, we get that $\cZ$ is the closure of the image by $\Phi$ of some open dense subset of $\HOM{A}(P_1, P_0)_{max}$.  Now, if $P_1 = Q \oplus Q_1$ and $P_0 = Q\oplus Q_0$, with $Q$ non-zero, then a generic morphism $f$ in $\HOM{A}(P_1, P_0)_{max}$ can be written as a direct sum $id_Q \oplus g$, for some $g: Q_1 \rightarrow Q_0$ (this follows from \cite[Theorem 5.2.2]{IOTW09}, see also \cite[Corollary 4.2]{DF09}).  Thus the generic morphism $f$ is not a minimal projective presentation of its cokernel; $g$ is the minimal presentation.  This is a contradiction; thus $P_1$ and $P_0$ have no direct factors in common.  This argument also shows that $P_1$ and $P_0$ determine $\cZ$; thus if $\cZ'$ is different from $\cZ$, the projectives $P'_1$ and $P'_0$ it determines are such that $(P_1, P_0) \ncong (P'_1, P'_0)$.
}

%-----------------------------------------------------------
\section{Generic value of cluster characters}\label{sect::generic value}
This section contains the main construction of the paper, namely that of a good candidate set for a basis of cluster algebras.  Before going into this construction, we recall background notions on cluster algebras and their categorification by triangulated categories.

%...............................................
\subsection{Cluster algebras and cluster categories}
\subsubsection{Cluster algebras}
Our main source for cluster algebras is the paper \cite{FZ07}.  An \emph{ice quiver} is a quiver $Q$ together with a set $F$ of vertices of $Q$ called \emph{frozen vertices}.  We will denote the non-frozen vertices by $1,2, \ldots, m$ and the frozen ones by $m+1, m+2, \ldots, n$.

Given an ice quiver, one can define a cluster algebra using the process of mutation.  Let $(Q,F)$ be an ice quiver without oriented cycles of length $1$ or $2$.  Let $i$ be a \emph{non-frozen} vertex.  Then the \emph{mutation of $(Q,F)$ at $i$} is the ice quiver $\mu_i(Q,F) = (Q',F')$ obtained from $(Q,F)$ in the following way:
\begin{enumerate}
	\item $F'=F$;
	\item for any subquiver $h\rightarrow i\rightarrow j$ of $Q$, add an arrow $h\rightarrow j$;
	\item reverse all arrows of $Q$ incident with $i$;
	\item delete all arrows of a maximal set of pairwise disjoint oriented cycles of length $2$.
\end{enumerate}

A \emph{seed} is a pair $\big(\underline{u}', (Q',F') \big)$, where $(Q',F')$ is an ice quiver without oriented cycles of length 1 or 2, and where $\underline{u}' = (u_1, \ldots, u_n)$ is an ordered subset of the field $\bQ(x_1, \ldots, x_n)$ which is free and generating.  The ordering $u_1, \ldots, u_n$ corresponds to the ordering of the vertices of $Q'$.

Let $i$ be an non-frozen vertex of $Q'$.  The \emph{mutation of the seed $\big(\underline{u}', (Q',F') \big)$ at $i$} is the seed $\mu_i\big(\underline{u}', (Q',F') \big) = \big(\underline{u}'', (Q'',F'') \big)$ defined by
\begin{itemize}
	\item $(Q'',F'') = \mu_i(Q',F')$;
	\item $u''_j = u'_j$ if $i \neq j$;
	\item $u''_i = (u'_i)^{-1}\Big(\prod_{j=1}^{n}(u'_j)^{(\# i\rightarrow j \textrm{ in } Q')} + \prod_{j=1}^{n}(u'_j)^{(\# j\rightarrow i \textrm{ in } Q')} \Big)$.
\end{itemize}

\begin{definition}[\cite{FZ02}]
Let $(Q,F)$ be an ice quiver as above.
\begin{itemize}
	\item The \emph{initial seed} associated to $(Q,F)$ is the seed $\big( (x_1, \ldots, x_n), (Q,F) \big)$.
	\item Any set $\underline{u}$ contained in a seed obtained by iterated mutations (at non-frozen vertices) from the initial seed is called a \emph{cluster}.
	\item Any element of a cluster is called a \emph{cluster variable}.
	\item A product of cluster variables contained in the same cluster is called a \emph{cluster monomial}.
	\item The \emph{cluster algebra} $\cA_{Q,F}$ associated to $(Q,F)$ is the sub-$\bQ$-algebra of $\bQ(x_1, \ldots, x_n)$ generated by all cluster variables.
	\item The \emph{upper cluster algebra} $\cA^+_{Q,F}$ \cite{BFZ05} associated to $(Q,F)$ is the sub-$\bQ$-algebra of $\bQ(x_1, \ldots, x_n)$ consisting of elements which can be written as a Laurent polynomial in the cluster variables of any given cluster.
\end{itemize}
\end{definition}

\begin{remarks}
\begin{enumerate}
	\item If the ice quiver has no frozen vertices, then $\cA_{Q,F} = \cA_Q$ is a cluster algebra \emph{without coefficients}.
	\item The original definition of \cite{FZ02} is more general, as it allows coefficients from any semifield and associates a cluster algebra to any skew-symmetrizable matrix (the datum of a quiver without oriented cycles of length 1 or 2 being equivalent to that of a skew-symmetric matrix). 
\end{enumerate}
\end{remarks}

\subsubsection{Cluster category of a quiver with potential}
Let $Q$ be a finite quiver.  It may have frozen vertices; this will not matter here.

Following \cite{DWZ08}, define a \emph{potential} $W$ on $Q$ as a (possibly infinite) linear combination of oriented cycles in $Q$.  More precisely, a potential is an element of the space $\widehat{\bC Q}/C$, where $\widehat{\bC Q}$ is the completed path algebra of $Q$ over $\bC$ and $C$ is the closure of the commutator subspace $[\widehat{\bC Q}, \widehat{\bC Q}]$.  The pair $(Q,W)$ is a \emph{quiver with potential}.

Given a potential $W$ and an arrow $a$ of $Q$, define the \emph{cyclic derivative} of $W$ with respect to $a$ as
\begin{displaymath}
	\partial_a W = \sum_{W = uav} vu.
\end{displaymath}
The \emph{Jacobian algebra} $J(Q,W)$ is the quotient of the completed path algebra $\widehat{\bC Q}$ by the closure of the ideal generated by the cyclic derivatives $\partial_a W$, as $a$ ranges over all arrows of $Q$.  The quiver with potential $(Q,W)$ is \emph{Jacobi-finite} if its Jacobian algebra is finite-dimensional.

Given a quiver with potential $(Q,W)$, a construction due to V.~Ginzburg yields a differential graded algebra $\Gamma$, called the \emph{completed Ginzburg dg algebra}.  We will not recall the construction of $\Gamma$ here; the reader is referred to \cite{G06} or \cite{Amiot08}.  Let us mention that $\Gamma$ is concentrated in non-positive degrees and that $H^0\Gamma = J(Q,W)$.

The \emph{perfect derived category} $\perf \Gamma$ of $\Gamma$ is the smallest triangulated category of the derived category $\cD\Gamma$ which contains $\Gamma$ and is closed under taking direct summands.  The category $\cD_{fd}\Gamma$ is the full subcategory of $\cD$ whose objects are those with finite-dimensional total cohomology.  By \cite[Theorem 2.17]{K08}, the category $\cD_{fd}$ is contained in $\perf\Gamma$.

\begin{definition}[\cite{Amiot08}]
The \emph{(generalized) cluster category} of $(Q,W)$ is the triangulated quotient
\begin{displaymath}
	\cC_{Q,W} = \perf\Gamma / \cD_{fd}\Gamma.
\end{displaymath}
\end{definition}

We will be working with the case where the cluster category is \emph{$\Hom{}$-finite}, which means that for any objects $X$ and $Y$ of $\cC_{Q,W}$, the vector space $\Hom{\cC}(X,Y)$ is finite-dimensional.

\begin{theorem}[\cite{Amiot08}]
The cluster category $\cC_{Q,W}$ is $\Hom{}$-finite if, and only if, $(Q,W)$ is Jacobi-finite.  In that case:
\begin{itemize}
	\item $\cC_{Q,W}$ is \emph{$2$-Calabi--Yau}, that is, for any objects $X$ and $Y$ of $\cC_{Q,W}$, there is a bifunctorial isomorphism
	  \begin{displaymath}
	    \Hom{\cC}(X,Y) \cong D\Hom{\cC}(Y, \Sigma^2 X),
    \end{displaymath}
where $D$ denotes the standard duality of vector spaces;

  \item $\Gamma$ is a \emph{cluster-tilting object} in $\cC_{Q,W}$, that is, $\Hom{\cC}(\Gamma, \Sigma\Gamma) = 0$, and for any object $X$ of $\cC_{Q,W}$, the condition $\Hom{\cC}(\Gamma, \Sigma X) = 0$ implies that $X$ is a finite direct sum of direct summands of $\Gamma$.
\end{itemize}
\end{theorem}

\subsubsection{The cluster character}
Let $(Q,F)$ be a finite ice quiver without oriented cycles of length 1 or 2, and let $W$ be \emph{non-degenerate} (see \cite{DWZ08}) potential on $Q$.  Then the cluster category $\cC_{Q,W}$ is known to categorify the cluster algebra $\cA_{Q,F}$, in a sense described below.

We will first need the notion of \emph{index} of an object of the cluster category.

\begin{proposition}[\cite{KR07}]
For any object $X$ of $\cC_{Q,W}$, there exists a triangle
\begin{displaymath}
	T_1 \longrightarrow T_0 \longrightarrow X \longrightarrow \Sigma T_1,
\end{displaymath}
with $T_1$ and $T_0$ in $\add \Gamma$.
\end{proposition}

\begin{definition}[\cite{DK08}]
The \emph{index} of $X$ (with respect to $\Gamma$) is the element
\begin{displaymath}
	\ind{\Gamma}X = [T_0] - [T_1] \in K_0(\add \Gamma)
\end{displaymath}
of the Grothendieck group $K_0(\add \Gamma)$. 
\end{definition}

Note that, if $\Gamma = \Gamma_1 \oplus \ldots \oplus \Gamma_n$ is a decomposotion of $\Gamma$ into indecomposable objects, then the group $K_0(\add\Gamma)$ is isomorphic to $\bZ[\Gamma_1]\oplus \ldots \oplus \bZ[\Gamma_n]$.

\begin{definition}[\cite{CC06}, \cite{Palu08}]
The \emph{cluster character} is the map
\begin{eqnarray*}
  CC: Obj(\cC_{Q,W}) & \longrightarrow & \bQ(x_1, \ldots, x_n) \\
           M & \longmapsto & x^{\ind{\Gamma}M} \sum_{\be \in \bN^{Q_0}} \Big( \chi \big( \Gr{e} (\Hom{\cC}(\Sigma^{-1}\Gamma, M))    \big) \Big) \prod_{j=1}^{n}\hat{y_j}^{e_j},
\end{eqnarray*}
where
\begin{itemize}
	\item $\Gr{\be}(X)$ is the quiver Grassmannian of the module $X$, that is, the projective variety of all submodule of $X$ with dimension vector $\be$;
	\item $\chi$ is the Euler characteristic;
	\item $x^{\ind{\Gamma}M}$ is the monomial $x_1^{g_1}\cdots x_n^{g_n}$, where $\ind{\Gamma}M = g_1[\Gamma_1] + \ldots + g_n[\Gamma_n]$; and
	\item $\hat{y_j} = \prod_{i=1}^{n} x_i^{(\# i \rightarrow j \textrm{ in } Q) - (\# j \rightarrow i \textrm{ in } Q)}$.
\end{itemize}
\end{definition}

There is a notion of mutation of objects in the cluster category, developped for any quiver with potential in \cite{KY09}.  An \emph{indecomposable reachable object} of $\cC_{Q,W}$ is an object which can be obtained from some $\Gamma_i$ by a sequence of mutations (at non-frozen vertices).

\begin{theorem}[\cite{CC06}, \cite{CK08}, \cite{Palu08}]
The cluster character induces a surjection from the set of isomorphism classes of indecomposable reachable objects of $\cC_{Q,W}$ to the set of cluster variables of $\cA_{Q,F}$.  
\end{theorem}

Note that the cluster character also allows to recover all the cluster monomials of $\cA_{Q,F}$.

%...............................................
\subsection{The $E$-invariant of Derksen--Fei in the cluster category}
Consider any cluster-tilting object $T$ in $\cC_{Q,W}$.  Then the functor $F= \Hom{\cC}(T, ?)$ induces an equivalence of categories \cite{KR07}
\begin{displaymath}
	\cC_{Q,W}/(\Sigma T) \longrightarrow \MOD \End{\cC}(T).
\end{displaymath}
This induces an equivalence $\add T \longrightarrow \proj\End{\cC}(T)$.
\begin{definition}
Let $T_0'$, $T''_0$, $T'_1$ and $T''_1$ be objects of $\add T$, and let $f'\in\Hom{\cC}(T'_1, T'_0)$ and $f''\in \Hom{\cC}(T''_1, T''_0)$.  Define $E(f',f'')$ to be the vector space $\Hom{K^b(\add T)}(\Sigma^{-1}f', f'')$.
\end{definition}
\begin{proposition}\label{prop::Einv in C}
\begin{enumerate}
	\item The space $E(f',f'')$ is isomorphic to the space $E(Ff', Ff'')$ of Derksen--Fei (see Definition \ref{defi::DF Einv}).
	\item Let $X'$ and $X''$ be cones of $f'$ and $f''$, respectively.  Then $E(f',f'')$ is isomorphic to $(T)(\Sigma^{-1}X', X'')$, the subspace of $\Hom{\cC}(\Sigma^{-1}X',X'')$ of morphisms which factor through an object of $\add T$.  In particular, $E(f',f'')$ only depends on the cones of $f'$ and $f''$.
\end{enumerate}
\end{proposition}
\demo{ Statement (1) is a direct consequence of the definition and of the equivalence $\add T \longrightarrow \proj\End{\cC}(T)$.  To prove statement (2), apply the functor $\Hom{\cC}(?, X'')$ to the triangle
\begin{displaymath}
	\xymatrix{ T'_1 \ar[r]^{f'} & T'_0 \ar[r] & X' \ar[r]& \Sigma T'_1
	}
\end{displaymath}
to get a diagram with exact rows and commuting left-most square
\begin{displaymath}
	\xymatrix{(T'_0, X'')\ar[r]^{(f')^*}\ar@{=}[d] & (T'_1, X'')\ar[r]^{\alpha}\ar@{=}[d] & (\Sigma^{-1}X', X'')\ar[r]^{\beta} & (\Sigma^{-1}T'_0, X'') \\
	 \Hom{B}(FT'_0, FX'')\ar[r]^{(Ff')^*} & \Hom{B}(FT'_1, FX'') \ar[r] & \Coker (Ff')^* \ar[r] & 0,
	}
\end{displaymath}
where $B = \End{\cC}(T)$ and we write $(W,Z)$ instead of $\Hom{\cC}(W,Z)$.  We know from \cite[Lemma 3.4]{DF09} that $E(Ff',Ff'')$ is isomorphic to $\Coker (Ff')^*$, which is in turn isomorphic to $\Coker (f')^*$.  This last space is isomorphic to $\Ima \alpha = \Ker \beta$.  But $\Ker \beta$ is exactly $(T)(\Sigma^{-1}X', X'')$.
}
\begin{remark}
This proposition enables us to formulate the theorem of Derksen--Fei (see Theorem \ref{theo::decomposition}) for any $\delta$ in $K_0(\add T)$ instead of in $K_0(\proj\End{\cC}(T))$.
\end{remark}

%...............................................
\subsection{Calabi--Yau reduction}\label{sect::reduction}
We recall properties of Calabi--Yau reduction in the sense of Iyama--Yoshino \cite{IY08}.  Let $(Q,F)$ be an ice quiver equipped with a non-degenerate potential $W$ such that $(Q,W)$ is Jacobi-finite.  We will write $\Gamma = \Gamma_N \oplus \Gamma_F$, where $\Gamma_N = \Gamma_1\oplus \ldots \oplus \Gamma_m$ is the direct sum of the indecomposable summands of $\Gamma$ corresponding to the non-frozen vertices, and $\Gamma_F = \Gamma_{m+1}\oplus\ldots\oplus \Gamma_n$ is the sum of those corresponding to the frozen vertices.  

Let $(\overline{Q}, \overline{W})$ be the quiver with potential without frozen vertices obtained from $(Q,W)$ in this way:
\begin{itemize}
	\item $\overline{Q}$ is the quiver obtained by removing all frozen vertices of $Q$ and all arrows incident to them;
	\item $\overline{W}$ is the potential obtained by removing all terms of $W$ involving arrows incident to frozen vertices.
\end{itemize}
We will denote by $\overline{\Gamma}$ the Ginzburg dg algebra of $(\overline{Q}, \overline{W})$.  Let $\cU$ be the full subcategory of $\cC_{Q,W}$ whose objects are those $X$ such that $\Hom{\cC}(X, \Sigma \Gamma_F) = 0$, and let $(\Gamma_F)$ be the ideal of all morphisms of $\cC_{Q,W}$ factoring through an object of $\add\Gamma_F$.

\begin{theorem}[\cite{IY08}]
The quotient $\cU/(\Gamma_F)$ is naturally endowed with a structure of triangulated category.
\end{theorem}

\begin{theorem}[Theorem 7.4 of \cite{K09}]\label{theo::keller reduction}
There is an equivalence of categories $\varphi : \cU/(\Gamma_F) \rightarrow \cC_{\overline{Q}, \overline{W}}$.  Under this equivalence, any $\Gamma_i$ with $i$ non-frozen is sent to $\overline{\Gamma}_i$.
\end{theorem}

Moreover, we have a way of comparing indices in $\cU \subset \cC_{Q,W}$ and in $\cC_{\overline{Q}, \overline{W}}$.  There is a natural surjection
\begin{displaymath}
	p : K_0(\add\Gamma) \longrightarrow K_0(\add\overline{\Gamma})
\end{displaymath}
sending $[\Gamma_i]$ to $[\overline{\Gamma}_i]$ if $i$ is non-frozen, and to $0$ otherwise.

\begin{proposition}\label{prop::indexcompare}
Let $X$ be an object of $\cU$.  Then its image $\varphi(X)$ in $\cC_{\overline{Q}, \overline{W}}$ under the equivalence of Theorem \ref{theo::keller reduction} has index $p(\ind{\Gamma}X)$.
\end{proposition}
\demo{ Recall from \cite{IY08} that triangles in $\cU/(\Gamma_F)$ are obtained in the following way.  Let $E\rightarrow X$ be a right $(\add\Gamma_F)$-approximation of $X$.  It embeds in a triangle $X\{-1\} \rightarrow E \rightarrow X \rightarrow \Sigma (X\{-1\})$.  The morphism $T_0 \rightarrow X$ yields a commutative diagram, where the rows are triangles:
\begin{displaymath}
	\xymatrix{ X\{-1\} \ar[r]\ar@{=}[d] & Q \ar[r]\ar[d] & T_0 \ar[r]\ar[d] & \Sigma (X\{-1\}) \ar@{=}[d] \\
	           X\{-1\} \ar[r] & E \ar[r] & X \ar[r] & \Sigma (X\{-1\}).
	}
\end{displaymath}
We will call this diagram $(*)$.  Then the sequence of morphisms $X\{-1\} \rightarrow Q \rightarrow T_0 \rightarrow X$ is sent to a triangle in $\cU/(\Gamma_F)$.

Applying the functor $\Hom{\cC}(\Gamma, ?)$ to $(*)$, we get a commutative diagram where rows and the left-most column are exact:
\begin{displaymath}
	\xymatrix{ \Hom{\cC}(\Gamma, T_0) \ar[r]^{u\phantom{xxxx}}\ar[d] & \Hom{\cC}(\Gamma,\Sigma(X\{-1\}))\ar[r]\ar@{=}[d] & \Hom{\cC}(\Gamma,\Sigma Q) \ar[r]\ar[d] & 0 \\
	           \Hom{\cC}(\Gamma, X)\ar[r]\ar[d] & \Hom{\cC}(\Gamma,\Sigma(X\{-1\}))\ar[r] & 0 & \\
	           0.
	}
\end{displaymath}
The map denoted by $u$, being a composition of two surjective maps, is itself surjective.  This implies that $\Hom{\cC}(\Gamma, \Sigma Q)$ vanishes.  Since $\Gamma$ is cluster-tilting, we get that $Q$ is in $\add\Gamma$.  Thus the index of $\varphi(X)$ is $[\varphi(T_0)] - [\varphi(Q)]$.  Moreover, the triangles in $(*)$ and \cite[Proposition 2.2]{Palu08} imply that $\ind{\Gamma}Q + \ind{\Gamma}\Sigma(X\{-1\}) = \ind{\Gamma}T_0$ and $\ind{\Gamma}E + \ind{\Gamma}\Sigma(X\{-1\}) = \ind{\Gamma}X$.  Thus we have $\ind{\Gamma}Q = \ind{\Gamma}T_0 - \ind{\Gamma}X + \ind{\Gamma}E = \ind{\Gamma}T_1 + \ind{\Gamma}E$.  Therefore $Q \cong T_1\oplus E$, and so the index of $\varphi(X)$ in $\cC_{\overline{Q}, \overline{W}}$ is $[\varphi(T_0)] - [\varphi(T_1\oplus E)] = p(\ind{\Gamma}X)$.
}

Now, given an object $Y$ of $\cC_{\overline{Q}, \overline{W}} \cong \cU/(\Gamma_F)$, we want to find an object $X$ of $\cU$ such that its image in $\cU/(\Gamma_F)$ is isomorphic to $Y$, and such that $X$ has no direct summands in $\add\Gamma_F$; we also want to know find an $(\add\Gamma)$-presentation of this $X$.

Let $\xymatrix{\overline{T_1} \ar[r]^{\overline{f}} & \overline{T_0}\ar[r] & Y\ar[r] & \Sigma\overline{T_1}  }$ be a triangle in $\cU/(\Gamma_F)$, with $\overline{T_1}$ and $\overline{T_0}$ in $\add\overline{\Gamma}$.  Let $T_1$ and $T_0$ in $\add\Gamma_N$ be lifts of $\overline{T_1}$ and $\overline{T_0}$, respectively, and let $f:T_1\rightarrow T_0$ be a lift of $\overline{f}$.  Note that a cone $Z$ of $f$ is not necessarily in $\cU$.  However, consider the following exact sequence:
\begin{displaymath}
	\xymatrix{\Hom{\cC}(T_0, \Gamma_F)\ar[r]^{f^*} & \Hom{\cC}(T_1, \Gamma_F)\ar[r] & \Hom{\cC}(\Sigma^{-1}Z, \Gamma_F)\ar[r] & 0.
	}
\end{displaymath}
Then the left $\End{\cC}(\Gamma_F)$- module $\Hom{\cC}(\Sigma^{-1}Z, \Gamma_F) \cong \Coker f^*$ has a projective cover, say $P_f$.  Let $T_f$ be an object of $\add\Gamma_F$ such that $P_f \cong \Hom{\cC}(T_f, \Gamma_F)$.

We have a commutative diagram of the form
\begin{displaymath}
	\xymatrix{ & \Hom{\cC}(T_f, \Gamma_F)\ar[d]\ar[dl]_{\alpha} & \\
	          \Hom{\cC}(T_1, \Gamma_F)\ar[r] & \Coker f^* \ar[r] & 0.
	}
\end{displaymath}  
The existence of the morphism $\alpha$ comes from the fact that $\Hom{\cC}(T_f, \Gamma_F)$ is projective.  Moreover, there exists an $h:T_1\rightarrow T_f$ such that $h^* = \alpha$.  Consider the map $(f,h)^t:T_1\rightarrow T_0\oplus T_f$.

\begin{lemma}\label{lemm::liftindex}
Let $X_f$ be a cone of $(f,h)^t:T_1\rightarrow T_0\oplus T_f$.
\begin{enumerate}
	\item $X_f$ lies in $\cU$ and has no direct summand in $\add\Gamma_F$.
	\item The image of $X_f$ in $\cU/(\Gamma_F)$ is isomorphic to $Y$.
	\item The isomorphism class of $X_f$ only depends on $\overline{f}$ (and not on the choice of the lift $f$).
	\item The isomorphism class of $T_f$ only depends on $\overline{f}$ (and not on the choice of the lift $f$).
\end{enumerate}
\end{lemma}
\demo{First, we prove (1).  We have an exact sequence
\begin{displaymath}
	\xymatrix{ \Hom{\cC}(T_0\oplus T_f, \Gamma_F)\ar[r]^{\phantom{xx}(f^*, h^*)} & \Hom{\cC}(T_1, \Gamma_F)\ar[r] & \Hom{\cC}(\Sigma^{-1}X_f, \Gamma_F)\ar[r] & 0.
	}
\end{displaymath}
By definition of $T_f$ and of $h$, the map $(f^*, h^*)$ is surjective.  Indeed, let $\beta\in \Hom{\cC}(T_1, \Gamma_F)$, and let $\overline{\beta}$ be its image in $\Coker f^*$.  Since $\Hom{\cC}(T_f, \Gamma_F)\rightarrow \Coker f^*$ is a projective cover, there exists a preimage $\gamma \in \Hom{\cC}(T_f, \Gamma_F)$ of $\overline{\beta}$, and $h^*(\gamma) = \alpha(\gamma) \in \Hom{\cC}(T_1, \Gamma_F)$ is also sent to $\overline{\beta}$ in $\Coker f^*$.  Thus $h^*(\gamma) - \beta \in \Ima f^*$, and so $\beta$ is in the image of $(f^*, h^*)$.  This proves the surjectivity of $(f^*, h^*)$, which implies that $\Hom{\cC}(\Sigma^{-1}X, \Gamma_F)$ vanishes, and this implies in turn that $X$ lies in $\cU$ by the definition of $\cU$.

Moreover, assume that $X_f \cong X'\oplus R$, with $R$ in $\add\Gamma_F$.  Then $R$ is a direct summand of $T_f$, since $X$ lies in $\cU$.  Write $T_f \cong Q\oplus R$.  We have a triangle
\begin{displaymath}
	\xymatrix{ T_1 \ar[r]^{(f, h_1, 0)^t\phantom{xxx}} & T_0\oplus Q\oplus R \ar[r] & X'\oplus R \ar[r] & \Sigma T_1,
	}
\end{displaymath}
which yields an exact sequence
\begin{displaymath}
	\xymatrix{ \Hom{\cC}(T_0\oplus Q\oplus R, \Gamma_F) \ar[r]^{(f^*, h_1^*, 0)} & \Hom{\cC}(T_1, \Gamma_F)\ar[r] & 0.
	}
\end{displaymath}
In particular, there is an epimorphism $\Hom{\cC}(Q, \Gamma_F)\rightarrow \Coker f^*$, and the first term is projective.  By minimality of the projective cover, we must have that $Q \cong T_f$, and so $R=0$.  Thus $X$ has no direct factor in $\add\Gamma_F$.

We now prove (2).  Let $E\rightarrow X$ be a right $(\add\Gamma_F)$-approximation of $X$.  We have a commutative diagram where rows are triangles:
\begin{displaymath}
	\xymatrix{ X\{-1\}\ar[r]\ar@{=}[d] & W \ar[r]\ar[d] & T_0\oplus T_f\ar[r]\ar[d] & \Sigma (X\{-1\})\ar@{=}[d] \\
	           X\{-1\}\ar[r]^{u} & E \ar[r] & X\ar[r] & \Sigma (X\{-1\}).
	}
\end{displaymath}
We know from the proof of Proposition \ref{prop::indexcompare} that $W \cong T_1 \oplus E$.  Moreover, the octahedral axiom of triangulated categories yields an octahedron

$\phantom{an octahedron}$\begin{xy} 0;<1pt,0pt>:<0pt,-1pt>:: 
(105,0) *+{X} ="0",
(74,90) *+{T_1} ="1",
(207,90) *+{T_0\oplus T_f} ="2",
(0,69) *+{E} ="3",
(103,147) *+{W} ="4",
(133,69) *+{\Sigma(X\{-1\})} ="5",
"0", {\ar|+"1"},
"2", {\ar"0"},
"3", {\ar"0"},
"0", {\ar@{.>}"5"},
"1", {\ar^{v}"2"},
"3", {\ar|+"1"},
"1", {\ar"4"},
"4", {\ar"2"},
"2", {\ar@{.>}"5"},
"4", {\ar^{0}"3"},
"5", {\ar@{.>}_{\Sigma u}|+"3"},
"5", {\ar@{.>}|+"4"},
\end{xy}

which implies that the morphism labelled by $v$ is isomorphic to $(f,h)^t$.  Therefore the map $W\rightarrow T_0\oplus T_f$ in the diagram above can be written as 
\begin{displaymath}
	\left( \begin{array}{cc}
f & *  \\
h & *  
\end{array} \right) : T_1\oplus E  \rightarrow  T_0 \oplus T_f.
\end{displaymath}
Thus its image in the quotient $\cU/(\Gamma_F)$ is $\overline{f}$.  Since the image of $X\{-1\} \rightarrow W\rightarrow T_0\oplus T_f \rightarrow X$ is a triangle in $\cU/(\Gamma_F)$, then the image of $X$ has to be $Y$.  This proves (2).

Statement (3) is a direct consequence of (1) and (2).  Statement (4) is a consequence of (3). 
}

%...............................................
\subsection{Generic values}
Let $(Q,F)$ be an ice quiver with no oriented cycles of length 1 or 2.  Let $W$ be a potential on $Q$ such that $(Q,W)$ is a Jacobi-finite non-degenerate quiver with potential.  Then C.~Amiot's cluster category $\cC_{Q,W}$ is $\Hom{}$-finite, $2$-Calabi--Yau and admits a cluster-tilting object $\Gamma$.  In this setting, the cluster character $CC$ is defined on the objects of $\cC_{Q,W}$.

For two objects $L$ and $M$ of $\cC_{Q,W}$, and for a morphism $\varepsilon$ from $L$ to $\Sigma M$, we denote by $mt(\varepsilon)$ any representative of the isomorphism class of ``middle terms'' $U$ in triangles
\begin{displaymath}
	M \longrightarrow U \longrightarrow L \longrightarrow \Sigma M.
\end{displaymath}
This notation is borrowed from \cite{Palu09}, as is the next result.  Recall (for instance, from sections 2.3 to 2.5 of \cite{EGA1}) that a \emph{locally closed subset} of a variety is the intersection of an open subset with a closed subset, that a \emph{constructible subset} of a variety is a finite union of locally closed subsets, and that a function from an algebraic variety to any abelian group is \emph{consctructible} if its image is finite and each fiber is a constructible subset of the variety.

\begin{proposition}[\cite{Palu09}]\label{prop::constr}
Let $L$ and $M$ be objects of $\cC_{Q,W}$.  Then the function
\begin{eqnarray*}
	\Hom{\cC}(L, \Sigma M) & \longrightarrow & \bQ(x_1, \ldots, x_n) \\
	      \varepsilon & \longmapsto & CC\big( mt(\varepsilon) \big)
\end{eqnarray*}
is constructible.
\end{proposition}
\demo{ This follows immediately from \cite[Proposition 9]{Palu09}.
}

Now, let $T_0$ and $T_1$ be objects in $\add \Gamma$.  It follows from Proposition \ref{prop::constr} that the function 
\begin{eqnarray*}
   \eta_{T_0, T_1} : \Hom{\cC}(T_1, T_0) & \longrightarrow & \bQ(x_1, \ldots, x_n) \\
      \varepsilon & \longmapsto & X'_{mt(\Sigma\varepsilon)}
\end{eqnarray*}
is constructible.  As in \cite{Dupont11}, we will define the map $I$ by using the fact that any constructible function on an irreducible variety admits a generic value (that is, there is a dense open subset of the domain of the function on which the function is constant).

Recall the decomposition $\Gamma = \Gamma_N \oplus \Gamma_F$, and our notation $\overline{\Gamma}$ for the image of $\Gamma$ in $\cU/(\Gamma_F)$.  Let $\overline{T_0}$ and $\overline{T_1}$ be objects of $\add\overline{\Gamma}$, $\overline{f}:\overline{T_1}\rightarrow \overline{T_0}$ be a morphism, $T_0$ and $T_1$ be lifts of $\overline{T_0}$ and $\overline{T_1}$ in $\add\Gamma_N$, and let $f:T_1\rightarrow T_0$ be a lift of $\overline{f}$.  Applying the functor $\Hom{\cC}(?, \Gamma_F)$, we get an exact sequence
\begin{displaymath}
	\xymatrix{\Hom{\cC}(T_0, \Gamma_F) \ar[r]^{f^*}  & \Hom{\cC}(T_1, \Gamma_F) \ar[r] & \Coker f^* \ar[r] & 0}
\end{displaymath}
of left $(\End{\cC}(\Gamma_F))$-modules.  Denote by $\overline{\delta}$ the element $[\overline{T_0}] - [\overline{T_1}]$ of $\add\overline{\Gamma}$.

\begin{lemma}\label{lemm::generic cover}
\begin{enumerate}
	\item There exists an object $T_F$ of $\add\Gamma_F$, which depends only on $\overline{T_1}$ and $\overline{T_0}$, such that for a generic morphism $\overline{f}:\overline{T_1}\rightarrow \overline{T_0}$, the left module $\Hom{\cC}(T_F, \Gamma_F)$ is a projective cover of $\Coker f^*$.
	\item Let $g$ be a generic morphism from $T_1$ to $T_0 \oplus T_F$, and let $X$ be its cone.  Then $\Hom{\cC}(X, \Sigma \Gamma_F) = \Coker g^* = 0$, and $X$ has no direct factor in $\add\Gamma_F$.
	\item Let $\delta = [T_0\oplus T_F] - [T_1] \in K_0(\add\Gamma)$.  If $T'_F$ lies in $\add\Gamma_F$, then $\delta + [T'_F] = \delta \oplus [T'_F]$ is a decomposition in the sense of Derksen--Fei (see Theorem \ref{theo::decomposition}).
\end{enumerate}
\end{lemma}
\demo{ We first prove (1).  There is a canonical surjection $\Hom{\cC}(T_1, T_0)\rightarrow \Hom{\overline{\cC}}(\overline{T_1}, \overline{T_0})$, through which a dense open subset of the first space is sent to something containing a dense open subspace of the first one.  Thus if $\overline{f}$ is generic, we can chose a lift $f$ which is also generic.  

Let $\Psi = \Hom{\cC}(?, \Gamma_F)$.  The image of the linear map $\Psi_{T_1, T_0}:\Hom{\cC}(T_1, T_0)\rightarrow \Hom{\End{}\Gamma_F}(\Psi(T_0), \Psi(T_1))$ is closed and irreducible.  Using the notations of section \ref{sect::cokernels}, there exists an $\bi$ such that $E_{\bi} \cap \Ima \Psi_{T_1, T_0}$ is open and dense in $\Ima \Psi_{T_1, T_0}$.  Let $\Phi : E_{\bi}\rightarrow \rep_{\bd}(\End{\cC}(\Gamma_F))$ be as in Lemma \ref{lemm::morphism}.  Then the composition $\Phi \circ \Psi_{T_1, T_0}$, when defined, sends any $f$ to $\Coker f^*$.  By Lemma \ref{lemm::CS}, the function
\begin{displaymath}
	\dim \Hom{\End{}\Gamma_F} (?, S) : \Ima \Phi \circ \Psi_{T_1, T_0} \longrightarrow \bN
\end{displaymath}
is upper semicontinuous for any simple module $S$.  Then it has a generic value, and an argument similar to that of Corollary \ref{coro::presentations} gives that all generic $f$ are such that $\Coker f^*$ have the same projective cover.  This projective cover is of the form $\End{\cC}(T_F, \Gamma_F)$, for some $T_F \in \add\Gamma_F$.  This proves (1).

Let $g$ be a generic morphism from $T_1$ to $T_0 \oplus T_F$, and let $X$ be a cone for $g$.  We have an exact sequence
\begin{displaymath}
	\xymatrix{ (T_0, \Gamma_F) \oplus (T_F, \Gamma_F) \ar[r]^{\phantom{xxxxxx}g^*} & (T_1, \Gamma_F) \ar[r] & (\Sigma^{-1}X, \Gamma_F) \ar[r] & 0,}
\end{displaymath} 
of left $(\End{\cC}\Gamma_F)$-modules, where we write $(Y,Z)$ instead of $\Hom{\cC}(Y,Z)$.  A generic $g$ minimizes the dimension of $(\Sigma^{-1}X, \Gamma_F)$; to prove (2), it is thus sufficient to show the existence of one $g$ such that $(\Sigma^{-1}X, \Gamma_F) \cong \Coker g^*$ vanishes.  Such a $g$ was constructed in Lemma \ref{lemm::liftindex} (1), where it was also shown that $X$ has no direct factor in $\add\Gamma_F$.  This proves (2).

To prove (3), let $\delta = [T_0\oplus T_F] - [T_1] \in K_0(\add\Gamma)$.  In view of the theorem of Derksen--Fei (see Theorem \ref{theo::decomposition}), it suffices to prove that $e(\delta, [T'_F]) = e([T'_F], \delta) = 0$.  Let $g\in \Hom{\cC}(T_1, T_0\oplus T_F)$ and $h\in \Hom{\cC}(0, T'_F)$ be generic morphisms.  Note that $h$ is zero.  By \cite[Lemma 3.4]{DF09}, $E(g,h)$ is isomorphic to the cokernel of
\begin{displaymath}
	\xymatrix{ \Hom{\cC}(T_0\oplus T_F, T'_F) \ar[r]^{g^*} & \Hom{\cC}(T_1, T'_F),
	}
\end{displaymath}
which is zero by part (2).  In the same way, we get that $E(h,g)$ is isomorphic to the cokernel of $h^*$, which vanishes since the target space of $h^*$ is zero.  This finishes the proof.
}

\begin{definition}\label{defi::I}
We define the map
\begin{displaymath}
	I : K_0(\add \overline{\Gamma}) \longrightarrow \bQ(x_1, \ldots, x_n)
\end{displaymath}
by letting $I([\overline{T_0}] - [\overline{T_1}])$ be the generic value of the map $\eta_{T_0\oplus T_F, T_1}$ defined above.
\end{definition}

In other words, $I([\overline{T_0}] - [\overline{T_1}])$ is the generic value taken by the cluster character $CC$ on cones of morphisms
\begin{displaymath}
	T_1 \longrightarrow T_0 \oplus T_F,
\end{displaymath}
where $T_F$ is defined as above.

\begin{proposition}\label{prop::linearly independant}
If the matrix of $Q$ is of full rank, then the elements in the image of $I$ are linearly independent over $\bZ$.
\end{proposition}
\demo{ Let $z_1, \ldots, z_r$ be distinct elements of the image of $I$, and let $p_1, \ldots, p_r$ be in $\bZ\bP$.  Suppose that $p_1z_1 + \ldots + p_rz_r = 0$.  Each $p_i$ is a Laurent polynomial in the frozen variables $x_{m+1}, \ldots, x_n$.  We can assume that the powers of the frozen variables in each term of each $p_i$ is non-negative; if that were not the case, we could multiply by positive powers of the frozen variables until it was so, without changing the equality to zero.

Next, we can also assume that each $p_i$ is a monomial in the frozen variables; if they were polynomials, we could use distributivity (the $z_i$ would then not be pairwise distinct; we will fix this below).

Finally, we can assume that each $p_i$ is actually an integer.  Indeed, let $p_i = \lambda x_{m+1}^{e_{m+1}} \cdots x_n^{e_n}$, and let $\overline{T_1}, \overline{T_0} \in \add(\overline{\Gamma})$ be such that $z_i = I([\overline{T_0}] - [\overline{T_1}])$.  Then $p_i z_i = \lambda x_{m+1}^{e_{m+1}} \cdots x_n^{e_n} I([\overline{T_0}] - [\overline{T_1}])$.  Using Lemma \ref{lemm::generic cover} (3), we get that this last term is equal to the cluster character $CC$ applied to the cone of a generic morphism 
\begin{displaymath}
	T_1 \longrightarrow T_0 \oplus T_F \oplus \bigoplus_{j=m+1}^{n}\Gamma_{j}^{e_{j}}.
\end{displaymath}

Note that, once that we have changed the $p_i$ to integers, we get a vanishing sum $\sum_{j=1}^{s} \lambda_j z'_j$, where the $\lambda_j$ are integers and the $z'_j$ are the values of the cluster character $CC$ applied to some objects.  Moreover, if $j \neq \ell$, then $z'_j$ and $z'_\ell$ are values of $CC$ applied to objects with distinct index.  This, allied to the fact that the matrix is of full rank, enables us to prove as in \cite[Corollary 4.4]{FK09} that the $\lambda_j$ must be zero.  Indeed, since the matrix is of full rank, there exists a grading on the ring of Laurent polynomial in $x_1, \ldots, x_n$ for which each $\hat{y}_i$ is of degree $1$.  Comparing lowest degrees of terms in the sum $\sum_{j=1}^{s} \lambda_j z'_j$ gives the vanishing of the $\lambda_j$.
}

%..........................................................
\subsection{Invariance under mutation}\label{sect::3}

We now define the mutation of indices.  This notion comes from the mutation of $Y$-variables of \cite{FZ07}, from the mutation of indices of \cite{DK08} and from the mutation of $\cX$-coordinates of \cite{FG09}.  Let $(Q,F)$ be an ice quiver, and $W$ be a potential on $Q$ such that $(Q,W)$ is non-degenerate and Jacobi-finite.

Let $\cX_{Q^{op}}(\bZ^t)$ be the set of tropical $\bZ$-points of the $\cX$-variety associated with the opposite quiver $Q^{op}$ (see \cite{FG09}).  Recall the isomorphisms 
\begin{displaymath}
	\cX_{Q^{op}}(\bZ^t) \cong \bZ^n \cong K_0(\add\Gamma).
\end{displaymath}
If $i$ is a non-frozen vertex of $Q$, then the mutation $\mu_i:\cX_{Q^{op}}(\bZ^t)\rightarrow \cX_{\mu_i(Q^{op})}(\bZ^t)$ is given by the tropicalization of the rule \cite[Formula (13)]{FG09}:
\begin{displaymath}
	Y'_j = \left\{ \begin{array}{ll}
             - Y_i & \textrm{if $i=j$;}\\
              Y_j - m[-Y_i]_+ & \textrm{if there are $m$ arrows from $i$ to $j$ in $Q$;}\\
              Y_j + m[Y_i]_+ & \textrm{if there are $m$ arrows from $j$ to $i$ in $Q$.}
\end{array} \right.
\end{displaymath}

On the other hand, there is a derived equivalence $\mu_i^-:\cC_{Q,W} \rightarrow \cC_{\mu_i(Q,W)}$ \cite{KY09} sending the object $\Gamma$ to the Ginzburg dg algebra $\Gamma'$ of $\mu_i(Q,W)$.  As seen in \cite[Section 4]{DK08} and the end of the proof of \cite[Proposition 2.7]{Plamondon09}, if an object $X$ has index $\ind{\Gamma}X = \sum_{j=1}^{n}g_j [\Gamma_j]$ in $\cC_{Q,W}$, then $\mu^-_i(X)$ has index $\ind{\Gamma'}X = \sum_{j=1}^{n}g'_j [\Gamma'_j]$ in $\cC_{\mu_i(Q,W)}$ given by
\begin{displaymath}
	g'_j = \left\{ \begin{array}{ll}
             - g_i & \textrm{if $i=j$;}\\
              g_j - m[-g_i]_+ & \textrm{if there are $m$ arrows from $i$ to $j$ in $Q$;}\\
              g_j + m[g_i]_+ & \textrm{if there are $m$ arrows from $j$ to $i$ in $Q$.}
\end{array} \right.
\end{displaymath}
We see from these formulas that the isomorphism $\cX_{Q^{op}}(\bZ^t) \cong K_0(\add\Gamma)$ commutes with mutation.

Recall that $(\overline{Q}, \overline{W})$ is obtained by removing all frozen vertices from $(Q,W)$, and that $\overline{\Gamma}$ is our notation for the Ginzburg dg algebra of $(\overline{Q}, \overline{W})$.  We have maps
\begin{displaymath}
	\xymatrix{ K_0(\add\Gamma) \ar@<1ex>[r]^p & K_0(\add\overline{\Gamma}) \ar@<1ex>[l]^j
	}
\end{displaymath}
where $p$ is a group homomorphism defined on the generators by $p([\Gamma_j]) = [\overline{\Gamma}_j]$ if $j$ is non-frozen and $0$ otherwise, and where $j([\overline{T}_0] - [\overline{T}_1]) = [T_0\oplus T_F] - [T_1]$, where $T_0$ and $T_1$ are lifts of $\overline{T}_0$ and $\overline{T}_1$ in $\add\Gamma_N$ and $T_F$ is as defined in Lemma \ref{lemm::generic cover}.

The next proposition will allow us to restrict ourselves to the case where the quiver $Q$ has no frozen vertices.

\begin{proposition}\label{prop::mutationreduction}
The following diagram commutes:
\begin{displaymath}
	\xymatrix{ \cX_{Q^{op}}(\bZ^t) \ar[r]^{\sim}\ar[d]_{\mu_i} & K_0(\add\Gamma)\ar[d]_{\mu_i} & K_0(\add\overline{\Gamma})\ar[l]_{j}\ar[d]_{\mu_i} \\
	           \cX_{\mu_i(Q^{op})}(\bZ^t) \ar[r]^{\sim} & K_0(\add\Gamma') & K_0(\add\overline{\Gamma'})\ar[l]_{j'}
	}
\end{displaymath}
\end{proposition}

We postpone the proof of Proposition \ref{prop::mutationreduction} until the end of the section.  Note that we have proved the commutativity of the left square in the discussion above.

Our main concern in this section is the proof of Theorem \ref{theo::3}.  We reformulate it thus:

\begin{theorem}[Reformulation of Theorem \ref{theo::3}]\label{theo::reform}
We have a commutative diagram
\begin{displaymath}
	 \xymatrix{ K_0(\add\overline{\Gamma})\ar[r]^{\mu_i}\ar[d]^{I} & K_0(\add\overline{\Gamma}')\ar[d]^{I} \\
	            \bQ(x_1, \ldots, x_n) \ar[r]^{\mu_i} & \bQ(x'_1, \ldots, x'_n),
	 }
\end{displaymath}
where the lower $\mu_i$ is the field isomorphism sending $x_j$ to $x'_j$ if $i\neq j$ and $x_i$ to $(x'_i)^{-1}\big( \prod_{\ell=1}^n (x'_\ell)^{\# i\rightarrow \ell} + \prod_{\ell=1}^n (x'_\ell)^{\# \ell\rightarrow i} \big)$.
\end{theorem}

\begin{remark}
In view of Proposition \ref{prop::mutationreduction}, it is sufficient to consider the case where $Q$ has no frozen vertices.  Indeed, we consider the diagram
\begin{displaymath}
	\xymatrix{ K_0(\add\overline{\Gamma}) \ar[d]^j\ar[r]^{\mu_i}\ar@/_4pc/[dd]_I & K_0(\add\overline{\Gamma}')\ar[d]^{j'}\ar@/^4pc/[dd]^I   \\
	           K_0(\add\Gamma) \ar[d]^I \ar[r]^{\mu_i} & K_0(\add\Gamma') \ar[d]^I \\
	           \bQ(x_1, \ldots, x_n) \ar[r]^{\mu_i} & \bQ(x_1', \ldots, x_n').
	}
\end{displaymath}
The left-most and right-most triangles commute by definition of $I$ and $j$.  The upper square commutes because of Proposition \ref{prop::mutationreduction}.  If Theorem \ref{theo::reform} was true for quivers without frozen vertices, then the lower square would commute, and so would the whole diagram, implying that Theorem \ref{theo::reform} is true for quivers with frozen vertices.
\end{remark}

We therefore restrict ourselves to the case where $Q$ has no frozen vertices in order to prove Theorem \ref{theo::reform}. 

We will need a lemma on the generic value of constructible functions.

\begin{lemma}\label{lemm::comm square constructible}
Let $W$, $X$ and $Y$ be irreducible varieties, and let $A$ be an abelian group.  Assume that we have a commutative diagram
\begin{displaymath}
	\xymatrix{ & W\ar[dr]^{v}\ar[dl]_{u}  & \\
	           X\ar[dr]_{\varphi} & & Y\ar[dl]^{\psi}  \\
	           & A &
	}
\end{displaymath}
where $u$ and $v$ are dominant morphisms of varieties and $\varphi$ and $\psi$ are constructible functions.  Then the generic values taken by $\varphi$ and $\psi$ are equal.
\end{lemma}
\demo{ Let $x\in A$ be the generic value taken by the function $\varphi$, and let $y\in A$ be that taken by $\psi$.  By definition, $\varphi^{-1}(x)$ is an open dense subset of $X$, and since $X$ is irreducible and $u$ is dominant, the intersection of $\varphi^{-1}(x)$ with the image of $u$ contains a dense open subset of $X$.  Thus $(\varphi \circ u)^{-1}(x)$ contains a dense open subset of $W$.  For similar reasons, $(\psi \circ v)^{-1}(y)$ contains a dense open subset of $W$.  Therefore $(\varphi \circ u)^{-1}(x)$ and $(\psi \circ v)^{-1}(y)$ have a non-empty intersection, and taking $w$ in their intersection, we get
\begin{displaymath}
	x = (\varphi \circ u)(w) = (\psi \circ v)(w) = y.
\end{displaymath}
This proves the result.
}

We can now prove Theorem \ref{theo::reform}

\demo{(of Theorem \ref{theo::reform}) Let $[T_0] - [T_1]$ be an element of $K_0(\add\Gamma)$.  We can assume that $T_0$ and $T_1$ have no direct factors in common.  Then $I([T_0] - [T_1])$ is, by definition, the generic value taken by the constructible function
\begin{eqnarray*}
	\eta_{T_0, T_1} : \Hom{\cC}(T_1, T_0) & \longrightarrow & \bQ(x_1, \ldots, x_n) \\
	\varepsilon & \longmapsto & CC(mt(\varepsilon)).	
\end{eqnarray*}

Let $\mu_i([T_0] - [T_1]) = [T'_0] - [T'_1]$ in $K_0(\add\Gamma')$.    Then $I([T'_0] - [T'_1])$ is the generic value taken by the constructible function
\begin{eqnarray*}
	\eta_{T'_0, T'_1} : \Hom{\cC'}(T'_1, T'_0) & \longrightarrow & \bQ(x'_1, \ldots, x'_n) \\
	\varepsilon & \longmapsto & CC(mt(\varepsilon)).	
\end{eqnarray*}

We want to show that $I([T'_0] - [T'_1])$ is the mutation of $I([T_0] - [T_1])$ at $i$; or, using our notation, that $I([T'_0] - [T'_1]) = \mu_i(I([T_0] - [T_1]))$.  It is a consequence of \cite[Corollary 4.14]{Plamondon10} that we have a commutative diagram
\begin{displaymath}
	\xymatrix{ \Hom{\cC}(\mu_i^+(T'_1), \mu_i^+(T'_0)) \ar[r]^{\phantom{xxxx}\mu_i^-}\ar[d]^{\eta} & \Hom{\cC'}(T'_1, T'_0)\ar[d]^{\eta} \\ 
	             \bQ(x_1, \ldots, x_n) \ar[r]^{\mu_i} & \bQ(x'_1, \ldots, x'_n),
	}
\end{displaymath}
where the two horizontal arrows are isomorphisms and where we omitted the indices of the maps $\eta$.  Thus $I([T'_0] - [T'_1]) = \mu_i(I([\mu_i^+(T'_0)] - [\mu_i^+(T'_1)]))$, and to prove the theorem it is therefore sufficient to show that $I([\mu_i^+(T'_0)] - [\mu_i^+(T'_1)]) = I([T_0] - [T_1])$.

We consider two cases.

\emph{Step 1: $\Gamma_i$ is not a direct summand of $T_1$.}  In that case, we can write $T_0 = \overline{T}_0 \oplus \Gamma_i^m$, where $\Gamma_i$ is not a direct summand of $\overline{T}_0$.  Recall that we have a (unique up to isomorphism) non-split triangle
\begin{displaymath}
	\xymatrix{\Gamma_i^* \ar[r]^{\alpha} & E' \ar[r]^{\beta} &  \Gamma_i \ar[r]^{\gamma} & \Sigma \Gamma_i^*
	}
\end{displaymath}
where $E' = \bigoplus_{a}\Gamma_{s(a)}$, the sum being taken over all arrows $a$ ending in $i$, and the morphism $\alpha$ is given by multiplication by these arrows on each coordinate.

  Then \cite[Proposition 2.7]{Plamondon09} (or rather, the triangle obtained at the end of its proof) allows us to write
\begin{displaymath}
  	\mu_i^+(T'_0) = \overline{T}_0 \oplus (E')^m \quad \textrm{and} \quad \mu_i^+(T'_1) = T_1 \oplus (\Gamma^*_i)^m.
\end{displaymath}

Consider the following diagram:
\begin{displaymath}
	\xymatrix{  \Aut{\cC}(\overline{T}_0 \oplus (E')^m) \times \Hom{\cC}(T_1, \overline{T}_0 \oplus (E')^m) \times \Aut{\cC}(T_1\oplus (\Gamma_i^*)^m)  \ar[d]_u \ar[dr]^{v}  &  \\
	           \Hom{\cC}(T_1 \oplus (\Gamma_i^*)^m, \overline{T}_0 \oplus (E')^m) \ar[dr]_{\eta}  & \Hom{\cC}(T_1, T_0) \ar[d]^{\eta}\\
	           & \bQ(x_1, \ldots, x_n).  
	}
\end{displaymath}
where $u$ and $v$ are defined as follows.  The morphism $v$ takes a triple $(g, f, g')$ and sends it to the composition $(id_{\overline{T}_0} \oplus \beta^{\oplus m})  \circ f$.  The morphism $u$ takes a triple $(g, f, g')$ and sends it to the morphism given in matrix form by
\begin{displaymath}
g  \left( \begin{array}{cc}
f_1 & 0  \\
f_2 & \alpha^{\oplus m}  
\end{array} \right) g' ,
\end{displaymath}
where $f =  \left( \begin{array}{c}
f_1   \\
f_2  
\end{array} \right)$.

If we could apply Lemma \ref{lemm::comm square constructible} to the above diagram, then the theorem would be proved.  Let us show that the hypotheses of the Lemma are fulfilled.  We easily see that the three varieties involved are irreducible, and we know from section \ref{sect::generic value} that the functions $\eta$ are constructible.  We must show that the square commutes and that $u$ and $v$ are dominant.

\emph{Substep 1 : the square commutes.} To show that the square commutes, we first notice that, for any $(g,f, g')\in \Aut{\cC}(\overline{T}_0 \oplus (E')^m) \times \Hom{\cC}(T_1, \overline{T}_0 \oplus (E')^m) \times \Aut{\cC}(T_1\oplus (\Gamma_i^*)^m)$, we have that $\eta v (g, f, g') = \eta v (id, f, id)$ (since $g$ and $g'$ do not occur in the definition of $v$) and that $\eta u(g, f, g') = \eta u  (id, f, id)$ (since the map $\eta$ takes the same value on orbits under the action of $\Aut{\cC}(\overline{T}_0 \oplus (E')^m) \times  \Aut{\cC}(T_1\oplus (\Gamma_i^*)^m)$).  Thus it is sufficient to show that $\eta v (id, f, id) = \eta u (id, f, id)$.  We invoke the octahedral axiom to get a diagram

$\phantom{an octahedron}$\begin{xy} 0;<1pt,0pt>:<0pt,-1pt>:: 
(105,0) *+{Y''\oplus \Sigma (\Gamma_i^*)^m} ="0",
(74,90) *+{T_1 \oplus (\Gamma_i^*)^m} ="1",
(207,90) *+{\overline{T}_0 \oplus \Gamma_i^m} ="2",
(0,69) *+{Y'} ="3",
(103,147) *+{\overline{T}_0 \oplus (E')^m} ="4",
(133,69) *+{\Sigma (\Gamma_i^*)^m} ="5",
"0", {\ar|+^{* \oplus id}"1"},
"2", {\ar_{}"0"},
"3", {\ar^{}"0"},
"0", {\ar@{.>}^{(a,b)}"5"},
"1", {\ar^{(v(id, f, id), 0)}"2"},
"3", {\ar|+_{}"1"},
"1", {\ar^{\!u(id, f, id)}"4"},
"4", {\ar_{id\oplus \beta^{\oplus m}}"2"},
"2", {\ar@{.>}_{\!\!\!(0, \gamma^{\oplus m})}"5"},
"4", {\ar^{}"3"},
"5", {\ar@{.>}_{}|+"3"},
"5", {\ar@{.>}^{}|+"4"},
\end{xy}

where $*$ is an unknown morphism, $Y'$ is the cone of $u(id, f, id)$ and $Y''$ is the cone of $v(id, f, id)$.  We need to show that $Y'$ and $Y''$ are isomorphic in order to show that the above square commutes, since $\eta u (id, f, id) = CC({Y'})$ and $\eta v (id, f, id) = CC(Y'')$.  The octahedron yields a commutative square
\begin{displaymath}
	\xymatrix{ Y''\oplus \Sigma (\Gamma_i^*)^m \ar[rr]^{(a,b)} \ar[d]^{*\oplus id} & & \Sigma(\Gamma_i^*)^m \ar[d]^{(0, -\Sigma\alpha^{\oplus m})^t} \\
	           \Sigma T_1 \oplus \Sigma (\Gamma_i^*)^m \ar[rr]^{-\Sigma u(id, f, id)} & & \Sigma \overline{T}_0 \oplus \Sigma (E')^m
	}
\end{displaymath}  
which, in turn, gives an equality of morphisms (in matrix form)
\begin{displaymath}
	\left( \begin{array}{cc}
0 & 0  \\
(-\Sigma \alpha^{\oplus m})\circ a & (-\Sigma\alpha^{\oplus m}) \circ b  
\end{array} \right)            = \left( \begin{array}{cc}
                                        * & 0  \\
                                        * & -\Sigma\alpha^{\oplus m}  
                                        \end{array} \right) ,
\end{displaymath}
where again the stars are unknown morphisms.  Thus $(-\Sigma\alpha^{\oplus m}) \circ b = -\Sigma\alpha^{\oplus m} $, and using the fact that the triangle $\xymatrix{\Gamma_i^* \ar[r]^{\alpha} & E' \ar[r]^{\beta} &  \Gamma_i \ar[r]^{\gamma} & \Sigma \Gamma_i^* }$ is a minimal $(\add\Gamma)$-copresentation of $\Gamma_i^*$, we get that $b$ is an isomorphism.  Therefore, the triangle (in the octahedron)
\begin{displaymath}
	\xymatrix{ Y'\ar[r] & Y''\oplus \Sigma (\Gamma_i^*)^m \ar[r]^{\phantom{xxx}(a,b)} & \Sigma (\Gamma_i^*)^m\ar[r] & \Sigma Y'
	}
\end{displaymath}
is isomorphic to a triangle
\begin{displaymath}
	\xymatrix{ Y'\ar[r] & Y''\oplus \Sigma (\Gamma_i^*)^m \ar[r]^{\phantom{xxx}(0,1)} & \Sigma (\Gamma_i^*)^m\ar[r] & \Sigma Y'
	}
\end{displaymath}
which is a direct sum of two triangles, of the form $\xymatrix{Y'\ar[r]& Y''\ar[r] & 0 \ar[r] & \Sigma Y'}$ and $\xymatrix{0\ar[r]& \Sigma (\Gamma_i^*)^m \ar[r] & \Sigma (\Gamma_i^*)^m \ar[r] & 0}$.  Thus $Y'$ and $Y''$ are isomorphic, and substep 1 is proven, that is, the above square commutes.

\emph{Substep 2 : the morphism $v$ is dominant}.  In fact, we show that $v$ is surjective, and thus dominant.  Indeed, let $f \in \Hom{\cC}(T_1, T_0)$.  Since $T_0 = \overline{T}_0 \oplus \Gamma_i^m$, we can write $f$ in matrix form as $f =  \left( \begin{array}{c}
f_1   \\
f_2  
\end{array} \right)$.

Now, since we have a triangle $\xymatrix{\Gamma_i^* \ar[r]^{\alpha} & E' \ar[r]^{\beta} &  \Gamma_i \ar[r]^{\gamma} & \Sigma \Gamma_i^* }$, and since the space $\Hom{\cC}(T_1, \Sigma\Gamma_i^*)$ vanishes (because $\Gamma_i$ is not a direct summand of $T_1$), we have that any morphism from $T_1$ to $\Gamma_i^m$ factors through $\beta^{\oplus m}$.  Thus we can write $f_2 = (\beta^{\oplus m})f'_2$, and we have a preimage of $f$ through $v$ of the form
\begin{displaymath}
	(id,  \left( \begin{array}{c}
f_1   \\
f'_2  
\end{array} \right) , id).
\end{displaymath}

\emph{Substep 3 : the morphism $u$ is dominant.} We will prove that the image of $u$ contains the following dense open subset of $\Hom{\cC}(T_1 \oplus (\Gamma_i^*)^m, \overline{T}_0 \oplus (E')^m)$:
\begin{displaymath}
	 \big\{  \left( \begin{array}{cc}
f_1 & h \\
f_2 & g\circ \alpha^{\oplus m} \circ g'
\end{array} \right)               \big| \ f_1, f_2, h \textrm{ are arbitrary}, \ g\in\Aut{\cC}((E')^m), \ g'\in \Aut{\cC}((\Gamma_i^*)^m)              \big\}.
\end{displaymath}
This subset is open because the subset $\{ g\circ \alpha^{\oplus m} \circ g' | \ g\in\Aut{\cC}((E')^m), \ g'\in \Aut{\cC}((\Gamma_i^*)^m) \}$ of $\Hom{\cC}((\Gamma_i^*)^m, (E')^m)$ is open, thanks to the fact that $\Gamma_i^m$ is rigid and to \cite[Lemma 2.1]{DK08}.  Let us show that it is contained in the image of $u$.  Let  
\begin{displaymath}
	\left( \begin{array}{cc}
	f_1 & h \\
f_2 & g\circ \alpha^{\oplus m} \circ g'
\end{array} \right) 
\end{displaymath}
be an element of it.  Then $h = h'\circ \alpha^{\oplus m}$ for some morphism $h'$, and we have
\begin{displaymath}
	\left( \begin{array}{cc}
	f_1 & h \\
f_2 & g\circ \alpha^{\oplus m} \circ g'
\end{array} \right)                    = \left( \begin{array}{cc}
	id & h' \\
0 & g
\end{array} \right)          \left( \begin{array}{cc}
	f_1 & 0 \\
g^{-1}f_2 & \alpha^{\oplus m} 
\end{array} \right)          \left( \begin{array}{cc}
	id & o \\
0 & g'
\end{array} \right),
\end{displaymath}
which is in the image of $u$.  Thus $u$ is dominant.

\emph{Substep 4}  We can now apply Lemma \ref{lemm::comm square constructible} to the above square, and as discussed earlier, this proves the theorem for the case considered in step 1.

\emph{Step 2 : $\Gamma_i$ is not a direct summand of $T_0$}.  In that case, we can write $T_1 = \overline{T}_1 \oplus \Gamma_i^n$, where $\Gamma_i$ is not a direct summand of $\overline{T}_1$.  We can use  arguments similar to those of step 1 to prove the theorem.  We could also work in the opposite triangulated category $\cC^{op}_{Q,W}$, and notice that
\begin{displaymath}
	\Hom{\cC}(T_1, T_0) = \Hom{\cC^{op}}(T_0, T_1),
\end{displaymath}
making step 2 in $\cC_{Q,W}$ equivalent to step 1 in $\cC^{op}_{Q,W}$.
}

We end this section by proving Proposition \ref{prop::mutationreduction}.  We will need the following result.
\begin{proposition}\label{prop::mutationdecomposition}
Let $\delta \in K_0(\add\Gamma)$, and let $\delta = \delta' \oplus \delta''$ be a decomposition in the sense of Derksen--Fei (see Theorem \ref{theo::decomposition}).  Then $\mu_i(\delta) = \mu_i(\delta')\oplus \mu_i(\delta'')$.
\end{proposition}
\demo{ Let $\delta' = [T'_0] - [T'_1]$ and $\delta'' = [T''_0] - [T''_1]$.  Consider the diagram
\begin{displaymath}
	\xymatrix{ & U' \times U''\ar[dr]^{v' \times v''}\ar[dl]_{u' \times u''} & \\
	           V' \times V'' \ar[dr]_{\dim E(?,?)} & & W' \times W'' \ar[dl]^{\dim E(?,?)} \\
	            & \bZ,  &
	}
\end{displaymath}
whose upper part is a product of two copies of the diagram in the beginning of the proof of Theorem \ref{theo::reform}.  More precisely,
\begin{itemize}
	\item $U' = \Aut{\cC}(\overline{T}'_0, (E')^{m'}) \times \Hom{\cC}(T'_1, \overline{T}'_0 \oplus (E')^{m'}) \times \Aut{\cC}(T'_1 \oplus (\Gamma_i^*)^{m'}) $;
	\item $U'' = \Aut{\cC}(\overline{T}''_0, (E')^{m''}) \times \Hom{\cC}(T''_1, \overline{T}''_0 \oplus (E')^{m''}) \times \Aut{\cC}(T''_1 \oplus (\Gamma_i^*)^{m''}) $;
	\item $V' = \Hom{\cC}(T'_1\oplus (\Gamma_i^*)^{m'}, \overline{T}'_0 \oplus (E')^{m'}  ) $;
	\item $V'' = \Hom{\cC}(T''_1\oplus (\Gamma_i^*)^{m''}, \overline{T}''_0 \oplus (E')^{m''}  ) $;
	\item $W' = \Hom{\cC}(T'_1, T'_0)  $;
	\item $W'' = \Hom{\cC}(T''_1, T''_0)  $; and
	\item $u'$, $u''$, $v'$ and $v''$ are defined as in the proof of Theorem \ref{theo::reform}. 
\end{itemize}
Let $(x',x'')$ be in $U'\times U''$.  It was shown in step 1.1 of the proof of Theorem \ref{theo::reform} that the cones of $u'(x')$ and $v'(x')$ are isomorphic, and that the cones of $u''(x'')$ and $v''(x'')$ are isomorphic.  Therefore, by Lemma \ref{prop::Einv in C} (2), we have that $\dim E(u'(x'), u''(x'')) = \dim E(v'(x'), v''(x''))$.  Thus the above square commutes.  Since we have already shown that $u'$, $u''$, $v'$ and $v''$ are dominant, and since we know from \cite[Corollary 3.8]{DF09} that $\dim E(?,?)$ is constructible, we can apply Lemma \ref{lemm::comm square constructible} to the above square: the generic values of the two occurences of $\dim E(?,?)$ coincide.  The one on the right is $e(\delta', \delta'')$, which vanishes, and the one on the left is $e(\mu_i(\delta'), \mu_i(\delta''))$, which is then zero.  Symmetric arguments show that $e(\mu_i(\delta''), \mu_i(\delta'))$ also vanishes.  Thus, by the theorem of Derksen--Fei (see Theorem \ref{theo::decomposition}), $\mu_i(\delta') \oplus \mu_i(\delta'')$ is a canonical decomposition.

There remains to be shown that $\mu_i(\delta)$ is the sum of $\mu_i(\delta')$ and $\mu_i(\delta'')$.  We know from Lemma \ref{lemm::signcoherence} that $\delta'$ and $\delta''$ are sign-coherent.  In that case, the definition of mutation of indices implies that $\mu_i(\delta' + \delta'') = \mu_i(\delta') + \mu_i(\delta'')$.  This finishes the proof.
}

\demo{(of Proposition \ref{prop::mutationreduction}).  There remains only to prove that the right-most square commutes.  Let $\overline{\delta}$ be in $K_0(\add\overline{\Gamma})$.  Lemma \ref{lemm::generic cover} gives us a characterization of $j'(\mu_i\overline{\delta})$: it is the element $z$ of $K_0(\add\Gamma')$ such that 
\begin{itemize}
	\item $p(z) = \mu_i(\overline{\delta})$;
	\item $z$ has no factor of the form $[R'_F]$, with $R'_F \in K_0(\add\Gamma'_F)$, in its canonical decomposition;
	\item if $R'_F \in K_0(\add\Gamma'_F)$, then $z + [R'_F] = z\oplus [R'_F]$.
\end{itemize}

  We see from the definition of mutation that $p'\circ \mu_i \circ j(\overline{\delta}) = \mu_i(\overline{\delta})$.  The fact that $\mu_i \circ j(\overline{\delta})$ has no term of the form $[R'_F]$ in its canonical decomposition, where $R'_F \in \add\Gamma'_F$, follows from Proposition \ref{prop::mutationdecomposition} and from the fact that $j(\overline{\delta})$ has no term of the form $[R_F]$ with $R_F \in \add\Gamma_F$ in its decomposition (see Lemma \ref{lemm::generic cover} (2)).  Now let $R'_F \in \add\Gamma'_F$.  Then $\mu_i(\mu_i \circ j(\overline{\delta}) + [R'_F]) = j(\overline{\delta}) + [R_F]$ since $j(\overline{\delta})$ and $[R_F]$ are sign-coherent, and $j(\overline{\delta}) + [R_F] = j(\overline{\delta}) \oplus [R_F]$ by Lemma \ref{lemm::generic cover} (3).  Thus $\mu_i \circ j(\overline{\delta}) + [R'_F] = \mu_i \circ j(\overline{\delta}) \oplus [R'_F]$ by Proposition \ref{prop::mutationdecomposition}.  This finishes the proof.
}

%----------------------------------------------------------------------
\section{Proof of Theorem \ref{theo::1}}\label{sect::1}
Consider the map $I: K_0(\add\overline{\Gamma}) \rightarrow \bQ(x_1, \ldots, x_n)$ defined in Definition \ref{defi::I}.  We proved in Proposition \ref{prop::linearly independant} that the elements in the image of $I$ are linearly independent over $\bZ\bP$.  

The fact that the image of $I$ is contained in the upper cluster algebra $\cA^+_Q$ follows from Theorem \ref{theo::3}.  Indeed, let $(u_1, \ldots, u_n)$ be a cluster obtained from the initial seed by a sequence of mutations at vertices $i_1, \ldots, i_s$.  Let $w$ be an element of $\bQ(x_1, \ldots, x_n)$, expressed in terms of the initial cluster.  Then its expression with respect to the cluster $(u_1, \ldots, u_n)$ is given by $\mu_{i_s}\circ \cdots \circ \mu_{i_1} (w)$.  Theorem \ref{theo::3} implies that, for any $\delta \in K_0(\add\Gamma)$, we have
\begin{displaymath}
	\mu_{i_s}\circ \cdots \circ \mu_{i_1} (I(\delta)) = I(\mu_{i_s} \circ \cdots \mu_{i_1}(\delta) ),
\end{displaymath}
which is a Laurent polynomial in the variables $u_1, \ldots, u_n$.  This being true for any cluster, $I(\delta)$ belongs to the upper cluster algebra.

Finally, suppose that $(Q,W)$ arises from the setting of \cite{GLS10}, and let $A$ be the Jacobian algebra of $(Q,W)$.  For any finite-dimensional representation $M$ of $A$, define $\psi_M$ thus: if
\begin{displaymath}
	\xymatrix{0\ar[r] & M\ar[r] & I_1\ar[r]^f & I_0}
\end{displaymath}
is a minimal injective presentation of $M$, then the injective modules $I_1$ and $I_0$ lift in the cluster category $\cC_{Q,W}$ to objects $\Sigma T_1$ and $\Sigma T_0$ of $\add \Sigma \Gamma$ through  $\Hom{\cC}(\Sigma^{-1}\Gamma, ?)$ (see \cite[Proposition 2.1]{KR07}).  Moreover, $f$ lifts to a morphism $\overline{f} \in \Hom{\cC}(\Sigma T_1, \Sigma T_0)$.  Then we put
\begin{displaymath}
	\psi_M = X'_{mt(\overline{f})}.
\end{displaymath}
For any irreducible component $\cZ$ of $\rep_{\bd}(A)$, let $\psi_{\cZ}$ be the generic value taken by $\psi_M$ in $\cZ$, and let
\begin{displaymath}
	null(\cZ) = \{ \bm \in \bN^{Q_0} \big| \ m_i = 0 \textrm{ if } d_i = 0 \}.
\end{displaymath}
Then it is proved in \cite[Theorem 5]{GLS10} that the set 
\begin{displaymath}
	B = \{ x^{\bm} \psi_{\cZ} \big| \ \cZ \textrm{ is strongly reduced in $\rep(A)$}, \bm \in null(\cZ) \}
\end{displaymath}
is a basis of the cluster algebra $\cA_Q$.  Let us prove that it is the image of the map $I$.

Assume that $\cZ$ is a strongly reduced component of $\rep_{\bd}(A)$.  Then, by Theorem \ref{theo::2}, we have $\cZ = \Psi([T_0] - [T_1])$ for some $[T_0]-[T_1] \in K_0(\add\Gamma)$.  By definition of $\Psi$ (see section \ref{sect::proof of 2}), $\cZ$ is the dual component of some $\cZ' = \Psi'([P_0]-[P_1])$, where $[P_0]-[P_1]\in K_0(\proj A^{op})$.  By definition of $\Psi'$, there is a dense open subset $\cU$ of $\Hom{A^{op}}(P_0, P_1)$ such that the union of orbits of cokernels of morphisms in $\cU$ contains a dense open subset of $\cZ'$.  Thus a generic representation in $\cZ'$ is isomorphic to a cokernel of a generic morphism in $\Hom{A^{op}}(P_0, P_1)$.  Dualizing, we get that a generic representation in $\cZ$ is isomorphic to a kernel of a generic morphism in $\Hom{A}(DP_1, DP_0)$. Note that $\Hom{A}(DP_1, DP_0)$ is isomorphic to $\Hom{\cC}(\Sigma T_1, \Sigma T_0)$ (since the $DP_i$ are finite-dimensional injective $A$-modules, see \cite[Proposition 2.1]{KR07}).  

Now, using Theorem \ref{theo::decomposition}, we get a canonical decomposition
\begin{displaymath}
	[T_0]-[T_1] = \delta_1 \oplus \ldots \oplus \delta_s,
\end{displaymath}
where the $\delta_i$ are indecomposable.  Assume that there are no non-negative terms in this decomposition.  This means that, generically in $\Hom{\cC}(\Sigma T_1, \Sigma T_0)$, $mt(\varepsilon)$ has no direct summand in $\add\Gamma$, so the generic value of $X'_{mt(\varepsilon)}$ is $\psi_{\cZ}$; in other words,
\begin{displaymath}
	I([T_0] - [T_1]) = \psi_{\cZ} \in B.
\end{displaymath}

Now, let $\bm \in null(\cZ)$.  Consider the non-negative element $[\bigoplus_{i=1}^{n}\Gamma_i^{m_i}] \in K_0(\add\Gamma)$. We will show that 
\begin{displaymath}
	[T_0]-[T_1] +[\bigoplus_{i=1}^{n}\Gamma_i^{m_i}] = \delta_1 \oplus \ldots \oplus \delta_s \oplus \bigoplus_{i=1}^{n}[\Gamma_i]^{\oplus m_i}
\end{displaymath}
is a canonical decomposition.  This will imply that 
\begin{eqnarray*}
I([T_0]-[T_1] +[\bigoplus_{i=1}^{n}\Gamma_i^{m_i}]) & = & I([\bigoplus_{i=1}^{n}\Gamma_i^{m_i}]) I([T_0]-[T_1]) \\
                                                    & = & x^{\bm}\psi_{\cZ},
\end{eqnarray*}
and will thus prove that the set $B$ is the image of the map $I$.

In order to do this, we work in the opposite category $\cC_{Q,W}^{op}$.  We use the functor $F = \Hom{\cC^{op}}((\Sigma^{op})^{-1}\Gamma, ?)$; in view of \cite[Theorem 4.4]{DF09}, we only need to show that for generic morphisms 
\begin{displaymath}
f'\in \Hom{A^{op}}(F(\Sigma^{op})^{-1}T_0, F(\Sigma^{op})^{-1}T_1) \textrm{ and } f''\in \Hom{A^{op}}(\bigoplus_{i=1}^{n}F(\Sigma^{op})^{-1}\Gamma_i^{m_i}, 0),
\end{displaymath}
the spaces $E(f',f'')$ and $E(f'',f')$ vanish.  Note that, by the above, a generic $f'$ has a cokernel $M''$ with dimension vector $\bd$.  Using \cite[Lemma 3.2]{DF09}, we thus get
\begin{eqnarray*}
 E(f',f'') & = & \Coker\big( \Hom{A^{op}}(F(\Sigma^{op})^{-1}T_1, 0)\rightarrow \Hom{A^{op}}(F(\Sigma^{op})^{-1}T_0, 0)  \big) \\
           & = & 0
\end{eqnarray*}
and
\begin{eqnarray*}
 E(f'',f') & = & \Coker\big( 0, M'')\rightarrow \Hom{A^{op}}(\bigoplus_{i=1}^{n}F(\Sigma^{op})^{-1}\Gamma_i^{m_i}, M'')  \big) \\
           & = & \Hom{A^{op}}(\bigoplus_{i=1}^{n}F(\Sigma^{op})^{-1}\Gamma_i^{m_i}, M'') \\
           & = & 0,
\end{eqnarray*}
since $\dimv M'' = \bd$, and $\bm \in null(\cZ)$.  This finishes the proof of the theorem.

\begin{remark}
Our proof that the image of the map $I$ is the set $B$ is valid for any Jacobi-finite quiver with potential, and not only for those arising from the setting of \cite{GLS10}.  Thus we can reformulate our results without mention of the cluster category by using representations of quivers with potential.

Let $(Q,W)$ be a Jacobi-finite quiver with potential, and let $A = J(Q,W)$.  Then we have a map
\begin{eqnarray*}
I: K_0(\add A) & \longrightarrow & B \\
  \delta & \longmapsto & x^{\bm}\psi_{\Psi(\delta)} 
\end{eqnarray*}
where $\bm \in \bN^{Q_0}$ is such that $\delta = \overline{\delta} \oplus [P_1^{m_1}] \oplus [P_n^{m_n}]$, and $\overline{\delta}$ has no non-negative summands in its canonical decomposition.  Then our results state that this map $I$ commutes with mutation, that $B$ lies in the upper cluster algebra, and that the elements of $B$ are linearly independent if the matrix of $Q$ is of full rank.  

Concretely, to compute $\Psi(\delta)$ when given an element $\delta = [P_1] - [P_0]$ of $K_0(\add A)$, one considers the corresponding injective $A$-modules $I_0$ and $I_1$. Then the kernels $M$ of generic morphisms $I_0 \rightarrow I_1$ form an open dense subset of $\Psi(\delta)$, and we have that $\psi_{\Psi(\delta)} = \psi_M$.
\end{remark}

%...............
\subsection{Link with a conjecture of V.~Fock and A.~Goncharov}
We now show how Theorems \ref{theo::1} and \ref{theo::3} are related to Conjecture 4.1 of \cite{FG09}.  Let $Q$ be a quiver without oriented cycles of length $\leq 2$, and let $\cA_Q$ be the associated cluster algebra (without coefficients).  The authors of \cite{FG09} conjecture, in a slightly different language, that there exists a bijection
\begin{displaymath}
	\cX_{Q^{op}}(\bZ^t) \longrightarrow E(\cA),
\end{displaymath}
where $E(\cA)$ is the subset of the cluster algebra consisting of elements which are Laurent polynomials with positive coefficients in the cluster variables of every cluster, and which cannot be written as a sum of two or more such elements.  This bijection should have the following properties:
\begin{enumerate}
	\item It should commute with mutation.
	
	\item A point $(a_1, \ldots, a_n)$ of $\cX_{Q^{op}}(\bZ^t)$ with non-negative coefficients should be sent to the element $\prod_{j=1}^n x_j^{a_j}$.
	
	\item The set $E(\cA)$ should be a $\bZ$-basis of the upper cluster algebra $\cA_Q^+$.
\end{enumerate}
Other conditions are described in \cite[Conjecture 4.1]{FG09}, but we will not discuss them here.
If we can equip the quiver $Q$ with a non-degenerate potential $W$ so that $(Q,W)$ is Jacobi-finite, then we have a good candidate for such a map.  

\begin{theorem}\label{theo::FG}
Let $(Q,W)$ be a non-degenerate, Jacobi-finite quiver with potential.  Then the map
\begin{displaymath}
	I : \cX_{Q^{op}}(\bZ^t) \cong \add\Gamma \longrightarrow \cA^+_Q
\end{displaymath}
defined in Definition \ref{defi::I} satisfies conditions 1 and 2 above.  If, moreover, $(Q,W)$ arises from the setting of \cite{GLS10}, then the image of $I$ satisfies condition 3.
\end{theorem}
\demo{ Condition 1 is Theorem \ref{theo::3}.  Condition 2 follows from the definition of $I$.  When we are in the setting of \cite{GLS10}, condition 3 follows from Theorem \ref{theo::1} and from the fact \cite[Theorem 3.3]{GLS10a} that in that case, the cluster algebra and the upper cluster algebra coincide.
}

Note that the coefficients of the elements in the image of $I$ need not be positive, as seen in \cite[Example 3.6]{DWZ09}.  Thus the image of $I$ is not contained in $E(\cA)$ in general.

The conjecture of \cite{FG09} discussed above is linked to one of \cite{GLS10}, where the authors conjecture that the set described in their Theorem 5 is a basis for the cluster algebra $\cA_Q$, starting from an arbitrary non-degenerate quiver with potential $(Q,W)$.  Using Theorem \ref{theo::2}, we know that, if $(Q,W)$ is Jacobi-finite, then this set of \cite{GLS10} is exactly the image of the map $I$, and by Example \ref{exam::conterexample} below, it is not necessarily contained in the cluster algebra, so that in the conjecture of \cite{GLS10}, one should replace ``cluster algebra'' by ``upper cluster algebra''.  If $I$ is a good candidate for the above map, then this is compatible with \cite[Conjecture 4.1]{FG09}.

\begin{example}\label{exam::conterexample}
The quiver with potential described below arises from the work of D.~Labardini-Fragoso \cite{Labardini09}\cite{Labardini11}.  We will show that the image of the map $I$ for this example is not contained in the associated cluster algebra, and that its cluster-category (which is $\Hom{}$-finite) has cluster-tilting objects which are not related by a finite sequence of mutations.

Consider the quiver
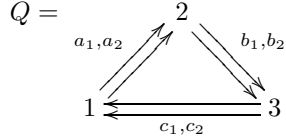
\begin{displaymath}
	Q = \xymatrix{ & 2\ar@<0.5ex>[dr]^{b_1,b_2}\ar@<-0.5ex>[dr] & \\
	              1\ar@<0.5ex>[ur]^{a_1,a_2}\ar@<-0.5ex>[ur] & & 3\ar@<0.5ex>[ll]^{c_1,c_2}\ar@<-0.5ex>[ll]
	}
\end{displaymath}
with potential $W= c_1b_1a_1 + c_2b_2a_2 - c_1b_2a_1c_2b_1a_2$.  As shown in \cite[Example 8.2]{Labardini11}, this quiver with potential is Jacobi-finite and non-degenerate.  Its Jacobian algebra $A$ is the path algebra of $Q$, modulo the relations
\begin{displaymath}
	\begin{array}{cccc}
	  c_1b_1 = c_2b_1a_2c_1b_2; & c_1b_1a_2 = 0; & b_1a_1c_2 = 0; & a_1c_1b_2 = 0; \\
	  c_2b_2 = c_1b_2a_1c_2b_1; & c_2b_2a_1 = 0; & b_2a_2c_1 = 0; & a_2c_2b_1 = 0; \\
	  b_1a_1 = b_2a_1c_2b_1a_2; & c_1b_2a_2 = 0; & b_1a_2c_2 = 0; & a_1c_2b_2 = 0; \\
	  b_2a_2 = b_1a_2c_1b_2a_1; & c_2b_1a_1 = 0; & b_2a_1c_1 = 0; & a_2c_1b_1 = 0; \\
	  a_1c_1 = a_2c_1b_2a_1c_2; & & & \\
	  a_2c_2 = a_1c_2b_1a_2c_1; & & & \\
  \end{array}
\end{displaymath}
and all non-alternating paths of length $4$ and all paths of length $7$ are zero.  These relations imply that $c_1b_1a_1 = c_2b_2a_2$; moreover, the Jacobian algebra is self-injective.

As a vector space, the indecomposable projective $A$-module $P_1 = e_1 A$ has a basis given by
\begin{displaymath}
	\{ e_1, c_1, c_2, c_1b_1, c_1b_2, c_2b_1, c_2b_2, c_1b_1a_1, c_1b_2a_1, c_2b_1a_2, c_1b_2a_1c_2, c_2b_1a_2c_1  \}.
\end{displaymath}
Similar calculations can be done for $P_2$ and $P_3$.  As a representation of the opposite quiver, we can draw $P_1$ as
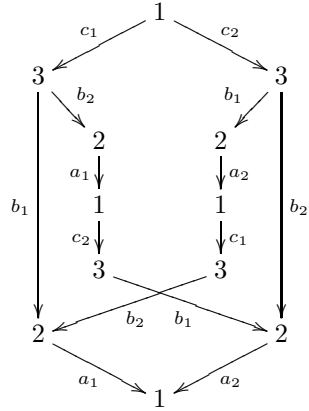
\begin{displaymath}
	\xymatrix@-1.0pc{ & & 1\ar[dll]_{c_1}\ar[drr]^{c_2} & & \\
	           3\ar[dddd]_{b_1}\ar[dr]^{b_2} & & & & 3\ar[dddd]^{b_2}\ar[dl]_{b_1} \\
	           & 2\ar[d]_{a_1} & & 2\ar[d]^{a_2} & \\
	           & 1\ar[d]_{c_2} & & 1\ar[d]^{c_1} & \\
	           & 3\ar[drrr]_{b_1} & & 3\ar[dlll]^{b_2} & \\
	           2\ar[drr]_{a_1} & & & & 2\ar[dll]^{a_2} \\
	           & & 1 & &	
	}
\end{displaymath}
Its socle is $S_1$; thus $P_1 = I_1$, the indecomposable injective at vertex $1$.  For similar reasons, $P_2 = I_2$ and $P_3 = I_3$.

Now, let $\cC$ be the cluster category of $(Q,W)$.  It is $\Hom{}$-finite, since $(Q,W)$ is Jacobi-finite (by results of \cite{Amiot08}).  Moreover, the functor $F = \Hom{\cC}(\Sigma^{-1}\Gamma, ?)$ sends $\Sigma^{-1}\Gamma_i$ to $P_i$ and $\Sigma\Gamma_i$ to $I_i = P_i$, for $i=1,2,3$.

Let us compute $I([\Gamma_1] - [\Gamma_3])$.  By definition, it is the generic value of the cluster character $X'_?$ applied to cones of morphisms in $\Hom{\cC}(\Gamma_3, \Gamma_1)$.  Equivalently, it is the value $X'_M$, where
\begin{displaymath}
	\xymatrix{ \Gamma_1\ar[r] & M\ar[r] & \Sigma\Gamma_3 \ar[r]^f & \Sigma \Gamma_1 }
\end{displaymath}
is a triangle and $f$ is generic in $\Hom{\cC}(\Sigma\Gamma_3, \Sigma\Gamma_1)$.  Applying the functor $F$, we get an injective presentation
\begin{displaymath}
	\xymatrix{ 0\ar[r] & FM\ar[r] & P_3 \ar[r]^{Ff} &  P_1, }
\end{displaymath}
where $Ff$ is generic.  Now, a generic $Ff$ in $\Hom{A}(P_3, P_1)$ is one for which $\dimv \Ker Ff$ is minimal.  We easily see that this minimal dimension vector is $(1,0,1)$ (for instance, one could take $Ff$ to be the left multiplication by $c_1 + c_2$), so that
\begin{displaymath}
	FM = \xymatrix{ & 0\ar@<0.5ex>[dl]\ar@<-0.5ex>[dl] & \\
	              \bC \ar@<0.5ex>[rr]^{\varphi_1,\varphi_2}\ar@<-0.5ex>[rr] & & \bC \ar@<0.5ex>[ul]\ar@<-0.5ex>[ul]
	}	
\end{displaymath}
with $\varphi_1\varphi_2 \neq 0$ as a representation of the opposite quiver.  Thus $FM$ has exactly $3$ subrepresentations, of dimension vectors $(0,0,0)$, $(0,0,1)$ and $(1,0,1)$.  Applying the cluster character, we get
\begin{eqnarray*}
 X'_M & = & x^{\ind{\Gamma}M} \sum_{e}\chi(\Gr{e}(FM)) \prod_{j=1}^{3}\hat{y}_j^{e_j} \\
      & = & \frac{x_1}{x_3} (1 + \hat{y}_3 + \hat{y}_1\hat{y}_3) \\
      & = & \frac{x_1}{x_3} (1 + x_1^{-2}x_2^{2} + x_2^{-2}x_3^2x_1^{-2}x_2^2) \\
      & = & \frac{ x_1^2 + x_2^2 + x_3^2}{x_1x_3}.
\end{eqnarray*}
This is the value of $I([\Gamma_1] - [\Gamma_3])$, and as shown in the proof of \cite[Proposition 1.26]{BFZ05}, it does not lie in the cluster algebra $\cA_Q$.

Now, the objects $\Gamma$ and $\Sigma\Gamma$ are cluster-tilting objects in $\cC$.  Note that 
\begin{displaymath}
	\ind{\Gamma}\Gamma = \sum_{j=1}^{3}[\Gamma_j] \quad \textrm{ and} \quad \ind{\Gamma}\Sigma\Gamma = - \sum_{j=1}^{3}[\Gamma_j].
\end{displaymath}
Let $Y$ be an object of $\cC$ with index $\sum_{j=1}^{3}y_j[\Gamma_j]$.  The index of $\mu^-_i(Y)$ is given by $\sum_{j=1}^{3}y'_j[\Gamma'_j] = \mu_i(\ind{\Gamma}Y) $, so that
\begin{displaymath}
	y'_j = \left\{ \begin{array}{ll}
             -y_i & \textrm{if $i=j$;}\\
              y_j + 2[y_i]_+ & \textrm{if there are  arrows from $j$ to $i$;}\\
              y_j - 2[-y_i]_+ & \textrm{if there are  arrows from $i$ to $j$.}
\end{array} \right.
\end{displaymath}
Thus we have that $y_i' = -y_i$, $y_{i+1}'= y_{i+1} - 2[-y_i]_+$ and $y_{i+2}'= y_{i+2} + 2[y_i]_+$, where $i,i+1,i+2$ are considered modulo $3$.  Thus
\begin{displaymath}
	y_i' + y_{i+1}' + y_{i+2}' = -y_i + y_{i+1} - 2[-y_i]_+ + y_{i+2} + 2[y_i]_+ = y_i + y_{i+1} + y_{i+2}.
\end{displaymath}
This shows that the sum of the coefficients appearing in the index of $Y$ is preserved under mutation of $Y$.  Since this sum is $3$ for $\Gamma$ and $-3$ for $\Sigma\Gamma$, the two objects cannot be related by a sequence of mutations.

\end{example}

\begin{remark}
In the above example, the fact that the sum of the coefficients of the indices is invariant under mutation was proved in \cite[Section 2.2]{FG11} (in a slightly different language).
\end{remark}

%----------------------------------------------------------------
\subsection*{Acknowledgements}

This work is an enhancement of the last chapter of my PhD thesis, supervised by Professor Bernhard Keller.  I wish to express my gratitude for his patience, his enthusiasm in discussing mathematics and his precious advice.  I would like to thank Professor Jan Schr\"oer for his valuable comments on a previous version of the paper.

%--------------------------------------------------------------------------
\bibliographystyle{amsplain} 
%\bibliography{cluster}

\begin{thebibliography}{00}

\bibitem{Amiot08}
Claire Amiot, \emph{Cluster categories for algebras of global dimension 2 and
  quivers with potential}, Ann. Inst. Fourier (Grenoble) \textbf{59} (2009),
  no.~6, 2525--2590.

\bibitem{ASS}
Ibrahim Assem, Daniel Simson, and Andrzej Skowro{\'n}ski, \emph{Elements of the
  representation theory of associative algebras. {V}ol. 1}, London Mathematical
  Society Student Texts, vol.~65, Cambridge University Press, Cambridge, 2006,
  Techniques of representation theory.

\bibitem{BFZ05}
Arkady Berenstein, Sergey Fomin, and Andrei Zelevinsky, \emph{Cluster algebras.
  {III}. {U}pper bounds and double {B}ruhat cells}, Duke Math. J. \textbf{126}
  (2005), no.~1, 1--52.

\bibitem{Borel}
Armand Borel, \emph{Linear algebraic groups}, second ed., Graduate Texts in
  Mathematics, vol. 126, Springer-Verlag, New York, 1991.

\bibitem{CC06}
Philippe Caldero and Fr{\'e}d{\'e}ric Chapoton, \emph{Cluster algebras as
  {H}all algebras of quiver representations}, Comment. Math. Helv. \textbf{81}
  (2006), no.~3, 595--616.

\bibitem{CK08}
Philippe Caldero and Bernhard Keller, \emph{From triangulated categories to
  cluster algebras}, Invent. Math. \textbf{172} (2008), no.~1, 169--211.

\bibitem{Cerulli09}
Giovanni Cerulli~Irelli, \emph{Cluster algebras of type ${A}_{\bf{2}}^{\bf
  {(1)}}$}, to appear in Algebras and Representation Theory, doi:
  10.1007/s10468-011-9275-5.

\bibitem{CE11}
Giovanni Cerulli~Irelli and Francesco Esposito, \emph{Geometry of quiver
  grassmannians of {K}ronecker type and applications to cluster algebras}, to
  appear in Algebra \& Number Theory, arXiv:1003.3037v2 [math.RT].

\bibitem{C93}
William Crawley-Boevey, \emph{Geometry of representations of algebras}, Lecture
  notes, 1993, available at the author's webpage.

\bibitem{CS02}
William Crawley-Boevey and Jan Schr{\"o}er, \emph{Irreducible components of
  varieties of modules}, J. Reine Angew. Math. \textbf{553} (2002), 201--220.

\bibitem{DK08}
Raika Dehy and Bernhard Keller, \emph{On the combinatorics of rigid objects in
  2-{C}alabi-{Y}au categories}, Int. Math. Res. Not. IMRN (2008), no.~11, Art.
  ID rnn029, 17.

\bibitem{DF09}
Harm Derksen and Jiarui Fei, \emph{General presentations of algebras},
  arXiv:0911.4913v1 [math.RA].

\bibitem{DWZ08}
Harm Derksen, Jerzy Weyman, and Andrei Zelevinsky, \emph{Quivers with
  potentials and their representations. {I}. {M}utations}, Selecta Math. (N.S.)
  \textbf{14} (2008), no.~1, 59--119.

\bibitem{DWZ09}
\bysame, \emph{Quivers with potentials and their representations {II}:
  applications to cluster algebras}, J. Amer. Math. Soc. \textbf{23} (2010),
  no.~3, 749--790.

\bibitem{DXX09}
Ming Ding, Jie Xiao, and Fan Xu, \emph{Integral bases of cluster algebras and
  representations of tame quivers}, to appear in Algebras and Representation
  Theory, doi: 10.1007/s10468-011-9317-z.

\bibitem{Dupont11}
Gr\'egoire Dupont, \emph{Generic cluster characters}, International Mathematics
  Research Notices IMRN (2011), doi: 10.1093/imrn/rnr024.

\bibitem{Dupont08}
\bysame, \emph{Generic variables in acyclic cluster algebras}, J. Pure Appl.
  Algebra \textbf{215} (2011), no.~4, 628--641.

\bibitem{DT11}
Gr\'egoire Dupont and Hugh Thomas, \emph{Atomic bases in cluster algebras of
  types ${A}$ and $\widetilde {A}$}, arXiv:1106.3758v1 [math.RA], 2011.

\bibitem{FG11}
Vladimir~V. Fock and Alexander~B. Goncharov, \emph{Cluster {X}-varieties at
  infinity}, arXiv:1104.0407v1 [math.AG].

\bibitem{FG09}
\bysame, \emph{Cluster ensembles, quantization and the dilogarithm}, Ann. Sci.
  \'Ec. Norm. Sup\'er. (4) \textbf{42} (2009), no.~6, 865--930.

\bibitem{FZ02}
Sergey Fomin and Andrei Zelevinsky, \emph{Cluster algebras. {I}.
  {F}oundations}, J. Amer. Math. Soc. \textbf{15} (2002), no.~2, 497--529
  (electronic).

\bibitem{FZ07}
\bysame, \emph{Cluster algebras. {IV}. {C}oefficients}, Compos. Math.
  \textbf{143} (2007), no.~1, 112--164.

\bibitem{FK09}
Changjian Fu and Bernhard Keller, \emph{On cluster algebras with coefficients
  and 2-{C}alabi-{Y}au categories}, Trans. Amer. Math. Soc. \textbf{362}
  (2010), no.~2, 859--895.

\bibitem{Gabriel74}
Peter Gabriel, \emph{Finite representation type is open}, Proceedings of the
  {I}nternational {C}onference on {R}epresentations of {A}lgebras ({C}arleton
  {U}niv., {O}ttawa, {O}nt., 1974), {P}aper {N}o. 10 (Ottawa, Ont.), Carleton
  Univ., 1974, pp.~23 pp. Carleton Math. Lecture Notes, No. 9.

\bibitem{GLS10a}
Christof Geiss, Bernard Leclerc, and Jan Schr{\"o}er, \emph{Kac--{M}oody groups
  and cluster algebras}, Advances in Mathematics \textbf{228} (2011), 329--433.

\bibitem{GLS10}
\bysame, \emph{Generic bases for cluster algebras and the {C}hamber {A}nsatz},
  J. Amer. Math. Soc. \textbf{25} (2012), 21--76.

\bibitem{G06}
Victor Ginzburg, \emph{Calabi--{Y}au algebras}, arXiv:math/0612139v3 [math.AG].

\bibitem{EGA1}
Alexandre Grothendieck and Jean Dieudonn\'e, \emph{{\'E}l\'ements de g\'eom\'etrie
  alg\'ebrique {I}}, Grundlehren der mathematischen Wissenschaften, vol. 166,
  Springer-Verlag, Berlin--Heidelberg, 1971.

\bibitem{HL10}
David Hernandez and Bernard Leclerc, \emph{Cluster algebras and quantum affine
  algebras}, Duke Math. J. \textbf{154} (2010), no.~2, 265--341.

\bibitem{IOTW09}
Kiyoshi Igusa, Kent Orr, Gordana Todorov and Jerzy Weyman, \emph{Cluster complexes via semi-invariants},
Compos. Math. \textbf{145} (2009), no. 4, 1001--1034. 

\bibitem{IY08}
Osamu Iyama and Yuji Yoshino, \emph{Mutation in triangulated categories and
  rigid {C}ohen-{M}acaulay modules}, Invent. Math. \textbf{172} (2008), no.~1,
  117--168.

\bibitem{Kashiwara90}
Masaki Kashiwara, \emph{Bases cristallines}, C. R. Acad. Sci. Paris S\'er. I
  Math. \textbf{311} (1990), no.~6, 277--280.

\bibitem{K09}
Bernhard Keller, \emph{Deformed {C}alabi--{Y}au completions}, to appear in
  Crelle's Journal, doi: 10.1515/CRELLE.2011.031.

\bibitem{K08}
\bysame, \emph{Calabi-{Y}au triangulated categories}, Trends in representation
  theory of algebras and related topics, EMS Ser. Congr. Rep., Eur. Math. Soc.,
  Z\"urich, 2008, pp.~467--489.

\bibitem{KR07}
Bernhard Keller and Idun Reiten, \emph{Cluster-tilted algebras are {G}orenstein
  and stably {C}alabi-{Y}au}, Adv. Math. \textbf{211} (2007), no.~1, 123--151.

\bibitem{KY09}
Bernhard Keller and Dong Yang, \emph{Derived equivalences from mutations of
  quivers with potential}, Advances in Mathematics \textbf{226} (2011),
  2118--2168.

\bibitem{Labardini11}
Daniel Labardini-Fragoso, \emph{Quivers with potentials associated to
  triangulated surfaces, {P}art {I}{I}: {A}rc representations},
  arXiv:0909.4100v2 [math.RT].

\bibitem{Labardini09}
\bysame, \emph{Quivers with potentials associated to triangulated surfaces},
  Proc. Lond. Math. Soc. (3) \textbf{98} (2009), no.~3, 797--839.

\bibitem{Lusztig90}
G.~Lusztig, \emph{Canonical bases arising from quantized enveloping algebras},
  J. Amer. Math. Soc. \textbf{3} (1990), no.~2, 447--498.

\bibitem{MSW11}
Gregg Musiker, Ralf Schiffler, and Lauren Williams, \emph{Bases for cluster
  algebras from surfaces}, arXiv:1110.4364v2 [math.RT], 2011.

\bibitem{Nak11}
Hiraku Nakajima, \emph{Quiver varieties and cluster algebras}, Kyoto J. Math.
  \textbf{51} (2011), no.~1, 71--126.

\bibitem{Palu09}
Yann Palu, \emph{Cluster characters {I}{I}: A multiplication formula}, to
  appear in Proc. LMS, arXiv:0903.3281v3 [math.RT].

\bibitem{Palu08}
\bysame, \emph{Cluster characters for 2-{C}alabi-{Y}au triangulated
  categories}, Ann. Inst. Fourier (Grenoble) \textbf{58} (2008), no.~6,
  2221--2248.

\bibitem{Plamondon10}
Pierre-Guy Plamondon, \emph{Cluster algebras via cluster categories with
  infinite-dimensional morphism spaces}, to appear in Compositio Mathematica, doi: 10.1112/S0010437X11005483.

\bibitem{Plamondon09}
\bysame, \emph{Cluster characters for cluster categories with
  infinite-dimensional morphism spaces}, Advances in Mathematics \textbf{227}
  (2011), no.~1, 1 -- 39.

\bibitem{SZ04}
Paul Sherman and Andrei Zelevinsky, \emph{Positivity and canonical bases in
  rank 2 cluster algebras of finite and affine types}, Mosc. Math. J.
  \textbf{4} (2004), no.~4, 947--974, 982.

\end{thebibliography}
\providecommand{\bysame}{\leavevmode\hbox to3em{\hrulefill}\thinspace}
%\providecommand{\MR}{\relax\ifhmode\unskip\space\fi MR }
%\providecommand{\MRhref}[2]{%
%  \href{http://www.ams.org/mathscinet-getitem?mr=#1}{#2}
%}
%\providecommand{\href}[2]{#2}

\end{document}